\documentclass[authoryear,revision,11pt]{elsarticle}
\usepackage{amssymb}
\usepackage{amsmath}
\usepackage{eurosym}
\usepackage{color}
\usepackage{bm}
\usepackage{booktabs}
\usepackage{multirow}
\usepackage{longtable}
\usepackage{arydshln}
\usepackage{geometry}
\usepackage{enumitem}
\usepackage{comment}

\usepackage{subcaption}
\usepackage{tikz}
\usepackage{pgfplots}
\usepackage{hyperref}

\usetikzlibrary{arrows}

\usepackage{graphicx}

\usepackage{etoolbox}
\makeatletter
\patchcmd{\ps@pprintTitle}{\footnotesize\itshape
Preprint submitted to \ifx\@journal\@empty Elsevier
\else\@journal\fi\hfill\today}{\relax}{}{}
\makeatother

\DeclareMathOperator*{\argmin}{arg\,min} % in your preamble 

\usepackage{amsthm,amsmath,eurosym}

\newtheorem{condition}{Condition}

\journal{European Journal or Operational Research}

\usepackage{lineno, blindtext}

\begin{document}

\begin{frontmatter}
    \title{A new scaling MCDA procedure putting together \\ pairwise comparison tables and the deck of cards method}

    \author[deb]{{\sc S. Corrente}}
    \author[ist]{{\sc J.R. Figueira}}
    %\author[deb,corl]{{\sc S. Greco}\corref{cor1}}\ead{salgreco@unict.it}
		\author[deb,corl]{{\sc S. Greco}}
    \address[deb]{Department of Economics and Business, University of Catania, Catania, Italy}
    \address[ist]{CEG-IST, Instituto Superior T\'{e}cnico,  Universidade de Lisboa, Portugal}
    \address[corl]{Portsmouth Business School, Centre of Operations Research and Logistics (CORL), \\
     University of Portsmouth,  Portsmouth, United Kingdom}

    \begin{abstract}
    \noindent This paper deals with an improved version of the deck of the cards method to render the construction of the ratio and interval scales more ``accurate''. The improvement comes from the fact that we can account for a richer and finer preference information provided by the decision-makers, which permits a more accurate  modelling of the strength of preference between two consecutive levels of a scale. Instead of considering only the number of blank cards between consecutive positions in the ranking of objects such as criteria and scale levels, we consider also the number of blank cards between non consecutive positions in the ranking. This information is collected in a pairwise comparison table that it is not necessarily built with precise values. We can consider also missing information and imprecise information provided in the form of ranges. Since the information provided by the decision-makers is not necessarily consistent, we propose also some procedures to help the decision-maker to make consistent his evaluations in a co-constructive way interacting with an analyst and reflecting and revising her/his judgments. The method is illustrated through and example in which, generalizing the SWING method, interacting criteria are aggregated through the Choquet integral.  
    \end{abstract}
    \vspace{0.25cm}
\end{frontmatter}

%\vfill\newpage
%\tableofcontents
%\vfill\newpage

\section{Introduction}\label{sec:Introduction}
\noindent Multiple Criteria Decision Aiding (MCDA) \citep{GrecoEhrgottFigueira2016} aims at supplying Decision-Maker (DM) with elements to reflect, conjecture, discuss, and argue about decisions in which a plurality of points of view are taken into consideration. MCDA procedures are based on an exchange of information between the DM, expressing relevant aspects of her/his preferences, and the analyst that uses such elements to build a decision model in a co-constructive approach involving the DM with her/his feedback in each stage. To construct the decision model, it is necessary to assign a value to several preference parameters depending on the nature of the adopted formal approach. Some examples of these parameters are the following:

\vfill\newpage

\begin{itemize}[label={--}]
    \item relative importance of criteria in outranking methods \citep{FigueiraEtAl2013};
    \item marginal value functions in multi-attribute value theory methods \citep{KeeneyRaiffa1976}; 
    \item capacities, that are interactive weights, permitting to represent interaction between criteria for the Choquet integral methods \citep{Choquet1953,Grabisch1996}.
\end{itemize}

An approach that has gained more and more attention for helping in eliciting parameters is the deck of card method (DCM) \citep{FigueiraRoy02}, which permits the DM to express differences in attractiveness by adding cards between consecutive elements. In this paper, we want to apply this method by using a richer information related to the difference of attractiveness between pairs of elements that are not necessarily consecutive. These differences in attractiveness are collected in a pairwise comparison table in which the element in line $r$ and column $s$ represents the number of cards corresponding to the difference in attractiveness between element $r$ and element $s$. 

Let us remember that the idea of pairwise comparison table has been largely used in MCDA in very well-known methods such as AHP \citep{Saaty1977} and MACBETH \citep{BanaCostaVansnick1994}. To the best of our knowledge pairwise comparison tables have been never coupled with a DCM. This paper proposes a methodology to fill this gap, using the pairwise comparison table to collect an information that is not necessarily complete, and which allows for some imprecision and some inconsistency. Moreover, since in case of imprecise and missing information, more than one comparison table can be compatible with the preference information provided by the DM, we propose to apply the Stochastic Multicriteria Accaptability Analsys (SMAA) \citep{LahdelmaHokkanenSalminen1998,PelissariEtAl2019,TervonenFigueira2008}. The application of SMAA in this case will permit to take into account all the comparison tables compatible with the information provided by the DM giving robust recommendations with respect to the problem at hand in probabilistic terms. 

After presenting the basic concepts, the main properties and the procedures to construct evaluations scales by means of pairwise comparison tables based on the DCM are presented. We show also how to apply this methodology to build a decision model in terms of the Choquet integral, that can be considered as a generalization of the usual weighted sum model, in case interaction of criteria has to be taken into consideration. In this regards we show that DCM can be seen as a support and refinement of the SWING method \citep{vonWinterfeldtEd86}. Finally, we illustrate how a decision process using pairwise comparison tables based on deck of cards method supports decisions by means of a didactic example.

The paper is organized as follows. Section 2 introduces the basics of the DCM approach; in Section 3 the DCM approach is used to build the capacity, that is, the weights of interactive criteria to be used to aggregate evaluations on single criteria through the Choquet integral; as a special case of our methodology we obtains a procedure to support the use of SWING methods to fix weights of criteria in the usual additive utility model.  The application of DCM approach  to assess interval scales is described in Section 4; Section 5 contains the description of the main concepts and definitions in our proposal; in Section 6, we explain how to restore consistency of the provided comparison table in case of inconsistent judgements; in Section 7, some extensions in the preference information provided by the DM in terms of imprecise and missing information are described; Section 8 contains a didactic example to which the proposed methodology is applied. Finally, some conclusions and further directions of research are gathered in the last section.

%%%%%%%%%%%%%%%%%%%%%%%%%%%%%%%%%%%%%%%%%%%%%%%%%%
\section{The basic DCM approach}\label{sec:Deck}%%%%
%%%%%%%%%%%%%%%%%%%%%%%%%%%%%%%%%%%%%%%%%%%%%%%%%%
\noindent In this section we describe the DCM approach used in SRF \citep{FigueiraRoy02} for the assessment of the weights of criteria for {\sc{Electre}} methods. Remember that \cite{FigueiraRoy02} proposed a modified version of Simos' DCM approach \cite[see][]{MaistreEtAl94} for determining the weights of criteria in {\sc{Electre}} methods and more generally for outranking based methods. For a list of applications of SRF see \cite{SiskosTsotsolas15}.  

In the remaining of this paper we will use the following basic notation. Let $A = \{a_1,\ldots,a_i,\ldots,a_m\}$ denote the set of alternatives to be assessed, and  $G = \{g_1,\ldots,g_j,\ldots,g_n\}$ denote the set of criteria. 

Criterion  $g_{j}\in G$ is a generic criterion to be maximized, while $a$ represents a generic alternative to be assessed on $g_{j}$; consequently, $g_{j}(a)$ is the performance of an alternative $a$ on a criterion $g_{j}$; moreover, $E_j = \{l_1,\ldots,l_k,\ldots,l_t\}$ is the scale of a discrete criterion $g_{j}$, where the scale levels are totally ordered, i.e., $l_1 \prec \cdots \prec l_k \prec \cdots \prec l_t$ ($\prec$ means ``strictly less preferred than'').

In addition, let  $U = \{u_1,\ldots,u_j,\ldots,u_n, \, u_{j}: E_j \rightarrow [0,1]\},$ denote the set of non-decreasing utility/value functions (one per criterion); $u_{j}$ denotes the generic utility/value function and $u_{j}(g_{j}(a))$ denotes the utility/value of the performance $g_{j}(a)$.

%%%%%%%%%%%%%%%%%%%%%%%%%%%%%%%%%%%%%%%%%%%%%%%%%%%
\subsection{Gathering the necessary preference information}\label{sec:Gathering_Ing}
%%%%%%%%%%%%%%%%%%%%%%%%%%%%%%%%%%%%%%%%%%%%%%%%%%%
\noindent In this subsection we show how the dialog between the analyst(s) and the decision-maker(s) or representatives(s) must be conducted to gather the necessary preference information. We will use an example with six criteria, $G = \{g_1, g_2, g_3, g_4, g_5, g_6\}$. The process of preference elicitation should be performed in a co-constructive interactive way between the analyst and the DM or a representative of the DM. 

    \begin{enumerate}
        \item in a first step, the analyst must prepare a set of cards representing considered criteria with their names on the corresponding cards  and, if necessary, some additional information (brief description, case study, criterion label, notation,...). An example of a criterion card is provided in Figure \ref{fig:criterion_card}. Then, the analyst should provide the DM with a first \emph{set of cards} and explain him the contents of each one.
        \begin{figure}[htb!]
            \centering
            \includegraphics[width=2.5cm, height=3.5cm]{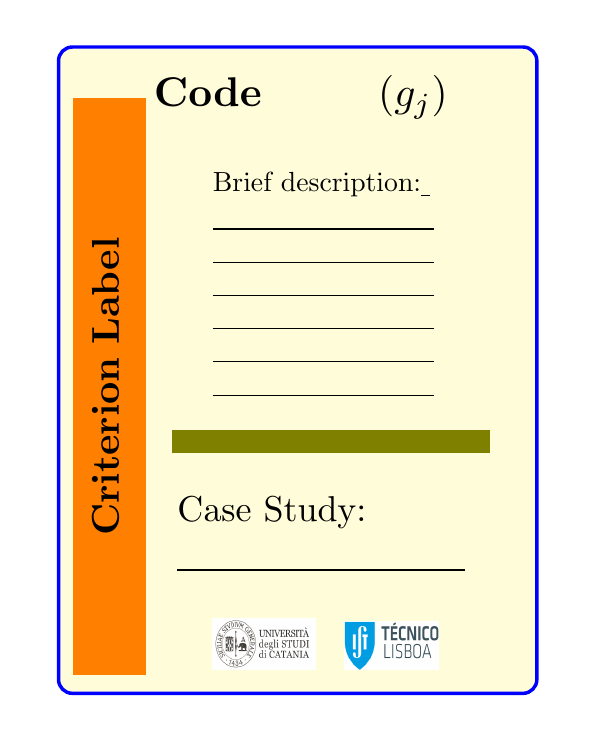}
            \caption{Example of a criterion card}
        \label{fig:criterion_card}
        \end{figure}
        \item in a second step, the analyst must also prepare a set of blank cards (as that one shown in Figure \ref{fig:Blank_card}) and provide them to the DM. The number of blank cards should be large enough, i.e., at least as much as the DM needs in the next step. 
         \begin{figure}[htb!]
            \centering
            \includegraphics[width=2.5cm, height=3.5cm]{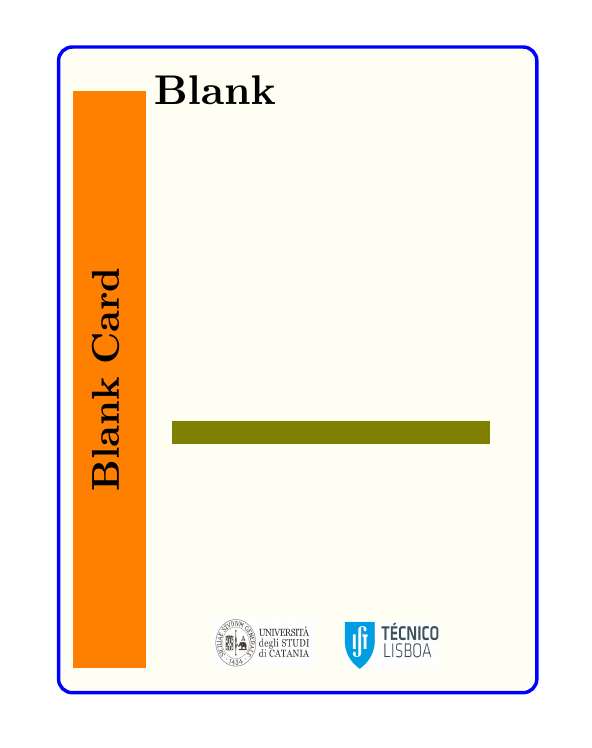}
            \caption{Example of a blank card}
        \label{fig:Blank_card}
        \end{figure}
        \item in a third step, the analyst must start to gather the preference information from the DM. The first question is to ask the DM to make a \emph{ranking} of the criteria cards from the least to the most important. Whenever the DM feels that some criteria are of the same importance he should put them in the same position in the ranking. Let us assume that in our case, the DM provided the following ranking:
        \[
            {\displaystyle \{g_3\} \prec \{g_4, g_5\} \prec \{g_1\} \prec \{g_2\} \prec \{g_6\}}.
        \]
        The least important criterion is $g_3$, then there are two with the same importance $\{g_4, g_5\}$, and so on, till the most important criterion, which is $g_6$. In our example, we have six criteria, but only five positions in the ranking since $g_4$ and $g_5$ are equally important.
        \item in the fourth step, an important information should be presented to the DM by the analyst, i.e., related to the fact that two consecutive positions in the ranking may be more or less close. In order to model this \textit{closeness} the DM can use blank cards and insert them in the consecutive intervals in the raking. No blank card does not mean that the criteria of the two consecutive positions have the same importance, but that this importance will be minimal (technically will represent the unit); one blank card means that the difference of importance is twice the minimal, and so on. Let us assume that the DM provided the following blank cards (in brackets) between the consecutive positions in the ranking.  
         \[
            {\displaystyle \{g_3\}\; [2]\; \{g_4, g_5\}\; [1]\; \{g_1\}\; [0]\; \{g_2\}\; [3]\; \{g_6\}}.
        \]
        
        Figure \ref{fig:ranking_blank} provides a more visual information about this step. 
     
        \begin{figure}[htpb!]
            \centering
            \includegraphics[width=16cm, height=5cm]{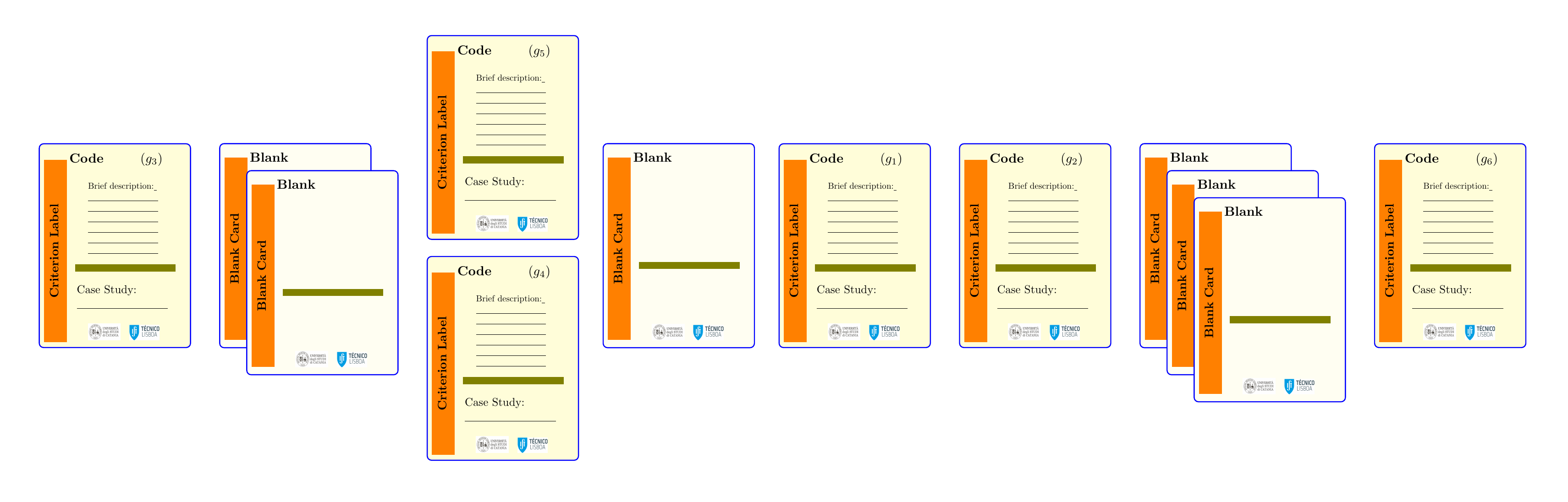}
            \caption{Ranking of criteria with blank cards}
            \label{fig:ranking_blank}
        \end{figure}

        \vfill\newpage
        
        \item finally, in the fifth step, the analyst should gather from the DM a crucial information for making possible the determination of the weights of criteria. Here, the DM must tell the analyst how many times the most important criterion (that one in the top position in the ranking) is more important that the least important one, i.e., that one in the lowest position in the ranking. In our example, this means how many times $g_6$ is more important than $g_3$. It is, in general, a difficult question. Sometimes, we can use the metaphor of vote and say that if criterion $g_3$ has one vote how many votes the DM could assign to criterion $g_6$? We do not have to work with only a precise value, sometimes two or three values or even a range can be provided and allow to have more sets of weights. Let us denote by $z$ (also called ratio $z$) this number. 
    \end{enumerate}

    Note that the information obtained in the fifth step is important in order to build a ratio scale. Indeed the nature of ratio scale characterizes the weights of {\sc{Electre}} methods \Citep{figueira2016electre} in which they have the meaning relative importance of criteria, while in MAVT (Multiple Attribute Value Theory ) based method \citep{KeeneyRaiffa1976} their meaning is of substitution rates of scaling factors.

%%%%%%%%%%%%%%%%%%%%%%%%%%%%%%%%%%%%%%%%%%%%%%%%%%%%%%%%%%%%%%%%
\section{The DCM approach for assessing a capacity for the Choquet integral}\label{sec:Ratio}%%%
%%%%%%%%%%%%%%%%%%%%%%%%%%%%%%%%%%%%%%%%%%%%%%%%%%%%%%%%%%%%%%%%
\noindent In many real world problems, the criteria are not mutually preferentially independent \citep{KeeneyRaiffa1976,Wakker1989} but present a certain degree of interaction. In general, we say that two criteria are positively (negatively) interacting if the importance assigned to them together is greater (lower) than the sum of the importance assigned to the two criteria singularly. To take into account this interaction, non-additive integrals are used in literature and, among them, the most well-known is the Choquet integral \citep{Choquet1953,Grabisch1996}. The Choquet integral is based on a capacity $\mu:2^G \rightarrow [0,1]$ being a set function that assigns a value not only to each criterion $g_i \in G$ but to all possible subsets of criteria $T \subseteq G$, so that monotonicity ($\mu(T_1)\leqslant\mu(T_2)$ for all $T_1\subseteq T_2\subseteq G$) and normalization ($\mu(\emptyset)=0$ and $\mu(G)=1$) constraints are satisfied. In this case, the formulation of the Choquet integral of the vector of utilities/values $\mathbf{u}(a)=(u_{1}(g_1(a)),\ldots,u_{n}(g_n(a))), a \in A$, is the following: 

\[
    {\displaystyle C_{\mu}(a) = \sum_{j=1}^{n}\Big(u_{(j)}(g_{(j)}(a)) - u_{(j-1)}(g_{(j-1)}(a))\Big)\mu({N_j}(a))}
\]
\noindent where $(\cdot)$ is a permutation of the indices $j, j=1, \ldots, n$,  such that $0=u_{(0)}(g_{(0)}(a))\leqslant u_{(1)}(g_{(1)}(a))\leqslant\ldots\leqslant u_{(n)}(g_{(n)}(a))$, and $N_{j}(a)=\left\{g_{i}\in G: u_{i}(g_{i}(a))\geqslant u_{(j)}(g_{(j)}(a))\right\}$. Let us observe that the use of differences $u_{(j)}(g_j(a)) - u_{(j-1)}(g_{(j-1)}(a))$ in the formulation of the Choquet integral imply that all evaluations $g_j(a)$, $g_j \in G$, are converted into the same scale through  $u_{1},\ldots,u_{n}$. In Section 4, generalizing what has been proposed in \citep{BotteroEtAl2013}, we shall explain how to assess utility functions $u_{1},\ldots,u_{n}$ using a DCM approach.  

To make things easier, in general, a M\"{o}bius transform of the capacity $\mu$ \citep{Rota1964} and $k$-additive capacities \citep{Grabisch1997} are used. In particular, the M\"{o}bius transform of $\mu$ is a set function $m$ such that $\mu(A)=\displaystyle\sum_{T\subseteq A}m(T)$, while a capacity $\mu$ is said $k$-additive iff its M\"{o}bius transform is such that $m(T)=0$ for all $T\subseteq G$ such that $|T|>k$. More details on this topic can be found in \cite{GrabischLabreuche2008a}. It has been observed that $2$-additive capacities are enough to represent preferences provided by the DM. Consequently, the Choquet integral of $a$ in case of $2$-additive capacities can be expressed in this linear form: 

\[
    {\displaystyle C_{m}(a) = \sum_{j=1}^{n}u_{j}(g_{j}(a))m\left(\left\{g_{j}\right\}\right)+\sum_{\left\{g_{i},g_{j}\right\}\subseteq G} \min\{u_{i}(g_{i}(a)),u_{j}(g_{j}(a))\}\cdot m\left(\left\{g_{i},g_{j}\right\}\right)}.
\]

As it is evident from the description above, to apply the Choquet integral one needs to assign values $\mu(T)$ to the subsets of criteria $T$ in $G$. For this reason, in the following, we shall describe the application of the DCM approach for assessing the capacities necessary to apply the Choquet integral following the procedure presented in \cite{BotteroEtAl18}. We start by considering the following dummy projects: 
    
    \begin{itemize}
        \item[$-$] projects $p_j$, for all $j=1,\ldots,n$, one per criterion, in which $p_j$ has the highest evaluation, say $u_j(g_j(p_j))=100$, on the criterion $g_j$ and the worst, say $0$, elsewhere, i.e., $\mathbf{u}(p_j) = (0,\ldots,0,u_j(g_j(p_j))=100,0,\ldots,0)$, for $j=1,\ldots,p$.    
        \item[$-$] $\vert O \vert$ dummy projects $p_{ij}$, related to the pairs $\{g_i,g_j\}$ contained in the set of pairs of interacting criteria $O\subseteq \{\{g_i,g_j\},g_i,g_j \in G\}$. The dummy projects $p_{ij}, \{g_i,g_j\} \in O$  have the maximal evaluations on the two interacting criteria $g_i$ and $g_j$, that is $u_i(g_i(p_{ij}))=u_j(g_j(p_{ij}))=100$, and the minimal evaluation on all the other criteria, that is $u_h(g_h(p_{ij}))=0$ for all $g_h\in G\setminus\{g_i,g_j\}$, so that $\mathbf{u}(p_{ij})=(0,\ldots,0,u_i(g_i(p_i))=100,0,\ldots,0,u_j(g_j(p_j))=100,0,\ldots,0)$, for all $\{i,j\} \in O$. The dummy projects $p_{ij}, \{g_i,g_j\} \in O$, are also denoted by $p_{k}$, $k=n+1,\ldots,n+\vert O \vert$, 
    \end{itemize}

    The method used to compute the capacities $\mu$ and, consequently, the M\"{o}bius transform $m$ of the capacity $\mu$, can be presented as follows: 

    \begin{enumerate}
        \item consider the reference projects: $P = \{p_1,p_2,\ldots,p_k,\ldots,p_t\}$ (where $t=n+\vert O \vert$);
        \item rank the projects in $P$ in different equivalent classes $R_1,\ldots,$ $R_h,\ldots,R_v$, from the worst to the best (that is, $R_1$ is the class containing the projects in the lowest position and $R_v$ the class containing the ones in the highest position). Let  $r_h$ denote a project representing the projects in $R_h$, for $h=1,\ldots,v$. Let $e_h$ denote the number of blank cards between $R_h$ and $R_{h+1}$, $h=1,\ldots,v-1$;
        \item assign a value to project $r_1$ (and consequently to all the projects in $R_1$), say $w(r_1)=\ell$  (it is frequent to consider $w(r_1)=1)$;
        \item insert blank cards $e_{h}$ between two successive subsets of projects $R_{h}$ and $R_{h+1}$ and provide the $z$ value representing how many times projects in $R_{v}$ are more important than projects in $R_{1}$;
        \item compute the value of each unit as follows:
        \[
            {\displaystyle \alpha = \frac{\ell(z-1)}{s}}
        \] where
        \[
            {\displaystyle s = \sum_{h=1}^{v-1}(e_h+1)};
        \]
        \item compute the values $w(r_h)$, $h=2,\ldots,v$, as follows:
        \[
            {\displaystyle w(r_h) = \ell+\alpha\left(\sum_{j=1}^{h-1}(e_h+1)\right)};
        \]
        \item compute the value of each project, $w(p_k) = w(r_h)$ for all $p_k \in R_h$, $h=1,\ldots,v$;
        \item compute the modified values $\overline{w}(p_k)$ as follows: $\overline{w}(p_k)=w(p_k)$ if $k=i\in G$;
        $\overline{w}(p_k)=w(p_k)-w(p_i)-w(p_j)$ if $p_k=p_{ij}, \{i,j\}\in O$, for $k \geqslant n+1$;
        \item compute the M\"{o}bius coefficients, $m_k$, and the capacities, $\mu_k$, for $k=1,\ldots,t$:
        \[
             {\displaystyle m_{k} = \frac{\overline{w}(p_{k})}{\displaystyle\sum_{j=1}^{t}\overline{w}(p_j)}},
        \]
        \noindent and
        \[
            {\displaystyle \mu_{k} = \frac{w(p_{k})}{\displaystyle\sum_{j=1}^{t}\overline{w}(p_j)}}.
        \]
        We have that $m_{k}=m\left(\{g_{k}\}\right)$ for $k=1,\ldots,n$ and $m_{k}=m\left(\{g_{i},g_{j}\}\right), \{g_{i},g_{j}\}\in O$ for $k\geqslant n+1$. Consequently, all $m_k$ must fulfill conditions $i^\prime)$ and $ii^\prime)$ below, ensuring monotonicity conditions and boundary conditions of capacity $\mu$, respectively: 
        \begin{description}
        \item[i)] for all $g_{j}\in G$ and for all $T\subseteq G\setminus\{g_{j}\}$, $m\left(\left\{g_{j}\right\}\right)+\displaystyle\sum_{g_{i}\in T, \{g_i,g_j\}\in O}m\left(\left\{g_{i},g_{j}\right\}\right)\geqslant 0$,
        \item[ii)] $m(\emptyset)=0$ and $\displaystyle\sum_{g_{j}\in G}m\left(\left\{g_{j}\right\}\right)+\sum_{\{g_{i},g_{j}\}\in O}m\left(\{g_{i},g_{j}\}\right)=1.$
        \end{description}
        Moreover, $m\left(\{g_{i},g_{j}\}\right)>0 \left[<0\right]$ if the DM stated that $g_{i}$ and $g_{j}$ are positively $\left[\mbox{negatively}\right]$ interacting. 
   \end{enumerate}
   
   Observe that
   \begin{itemize}
       \item in case there are not interacting criteria, the above procedure amounts to apply the DCM approach to the wellknown SWING method \citep{vonWinterfeldtEd86},
   \item consequently, the above procedure can be seen also as an extension of the SWING method to the case of aggregation of evaluations of interacting criteria through the Choquet integral.
     \end{itemize}

\section{The DCM approach for assessing utilities in interval scales}\label{sec:Deck}
\noindent  The utility values of the MAVT and Choquet integral methods are the levels of a common interval scale, in general, within the range $[0,1]$. The translation from the original scales of the criteria to a single common interval scale requires the use of a procedure that should account for the intensity of preferences between consecutive intervals of the scale. In this section, we present a new procedure for defining an interval scale based on the concepts of the deck of cards method \citep{FigueiraRoy02}. The procedure presented here allows scales not necessarily within the range $[0,1]$ to be constructed.

    In order to build an interval scale we need to define at least two reference levels (instead of the definition of $z$, as in the case of ratio scales), to anchor the computations. If more than two reference levels are defined, we can replicate the procedure for every two consecutive reference levels.

    \begin{enumerate}
        \item consider the scale $E_j = \{l_1,\ldots,l_k,\ldots,l_t\}$; 
        \item define two reference levels, $l_p$ and $l_q$ and define their utilities. It is frequent to use $u(l_p)=0$ and $u(l_q)=1$;
        \item insert the blank cards in the successive interval of the rankings: 
            \[
                {\displaystyle l_1\;e_1\;l_2\;\cdots\;l_p\;e_p\;l_{p+1}\;e_{p+1}\;\cdots\;l_{k}\;e_{k}\;l_{k+1}\;\cdots\;l_{q-1}\;e_{q-1}\;l_q
                \;\cdots\;l_{t-1}\;e_{t-1}\;l_t};
            \]
        \item consider only the levels in between $l_p$ and $l_q$  and determine the value of the unit:
            \begin{equation}\label{AlphaValue}
                {\displaystyle \alpha = \frac{u(l_q)-u(l_p)}{h}},
            \end{equation} 
            
            \noindent where
            
            \begin{equation}\label{BlankCards}
                {\displaystyle h = \sum_{r=p}^{q-1}(e_r+1)},
            \end{equation}
            
            \noindent which is the number of units between levels $l_p$ and $l_q$;
        \item compute the utility value, $u(l_k)$, for each level, $k=1,\ldots,t$, as follows:
            \[
            u(l_k) = \left\{
                \begin{array}{ccl}
                    {\displaystyle u(l_p)-\alpha\left(\sum_{j=k}^{p-1}(e_j+1)\right),} & \mbox{for} & k=1,\ldots,p-1,  \\
                    & & \\
                    {\displaystyle u(l_p)+\alpha\left(\sum_{j=p}^{k-1}(e_j+1)\right),} & \mbox{for} & k=p+1,\ldots,q-1,q+1,\ldots,t. \\
                \end{array}
                \right.
            \]
    \end{enumerate}

    An attempt to build an interval scale was also proposed by \cite{PictetBollinger08}, but without considering the value of the unit and considering instead the blank cards as positions in the ranking. This makes the method unpractical when confronted with a large number of blank cards in the intervals between consecutive levels. In addition, the formula to compute the values of each level cannot assign values strictly lower than the value of the lowest reference level.

%%%%%%%%%%%%%%%%%%%%%%%%%%%%%%%%%%%%%%%%%%%%%%%%%%%%%%%%%%%%%%%%%%%%%%%%%%%%%%%%%%%%%%%%%%%%%%
\section{Pairwise comparison tables based on the deck of cards method }\label{PCTDCM}%%
%%%%%%%%%%%%%%%%%%%%%%%%%%%%%%%%%%%%%%%%%%%%%%%%%%%%%%%%%%%%%%%%%%%%%%%%%%%%%%%%%%%%%%%%%%%%%%
\noindent The DCM approach, as described in the literature and recalled in the previous sections, is based on the comparisons between the elements in consecutive levels. Indeed, the DM is asked to rank order the cards corresponding to the considered elements and, after, to add blank cards between two consecutive cards to express the difference in the appreciation of elements in one level and elements in the following level. However, one can imagine that the DM could enrich the information he supplies by expressing the differences in the appreciation not only between consecutive levels of elements, but also between non consecutive levels. This information can be collected in a pairwise comparison table that has a nature similar to the pairwise comparison tables considered in two well known MCDA methods, AHP \citep{Saaty1977} and MACBETH \citep{BanaCostaVansnick1994}. The pairwise comparison tables of the DCM approach has the advantage of a visual support represented by the cards, that can aid the DM in defining and expressing his values and preferences. This seems an important point in the perspective of an MCDA methodology that has the primary scope of supporting the DM in discussing and arguing  in order to construct a conviction on the decision to be taken. In the following we shall introduce and discuss in deep the theory and the practice of the pairwise comparison table based on the DCM approach.

%%%%%%%%%%%%%%%%%%%%%%%%%%%%%%%%%%%%%%%%%%%%%%%%%%%%%%%%%%%%%%%%%%%%%%%%%%%%%%%%%%%%%%%%%%%%%
\subsection{Basic elements and consistency condition}\label{RestoringConsistency}%%%
%%%%%%%%%%%%%%%%%%%%%%%%%%%%%%%%%%%%%%%%%%%%%%%%%%%%%%%%%%%%%%%%%%%%%%%%%%%%%%%%%%%%%%%%%%%%%
\noindent Consider the scale levels as in the previous sections, $l_{1},\ldots,l_{t}$, where $l_{1}$ is the lowest level and $l_{t}$ is the highest one. In the classic DCM approach, the DM only introduces cards between consecutive levels. Now, we will enrich the preference information asking the DM to fulfill a comparison table, $C$, where each entry, $e_{pq}$, denotes the number of blank cards that should be inserted between levels $l_{p}$ and $l_q$. Of course, as in the classical DCM approach, the greater the number of blank cards between $l_{p}$ and $l_{q}$, the greater the strength or intensity of preference between the two considered levels. The comparison table has to fulfill an important consistency condition: 

\begin{condition}[Consistency]\label{Precies Consistency Condition} Given the comparison table, $C$,
    
    \[
\left\vert
\begin{array}{c|ccccccccc}\hline
         & l_{1} & \ldots & l_p & \ldots & l_k & \ldots & l_q & \ldots & l_t \\ \hline
 l_{1}  &       \blacksquare &  &  &  &  &  &  &  & \\
 \vdots & &      \blacksquare &  &  &  &  &  &  & \\
 l_{p}  & & &     \blacksquare &  & e_{pk} & \ldots & e_{pq} &  & \\
 \vdots & & & &    \blacksquare &  & \ddots & \vdots &        & \\
 l_{k}  & & & &  &  \blacksquare &   & e_{kq} &  & \\
 \vdots & & & & &  & \blacksquare  &  &  & \\
 l_{q}  & &  & && & & \blacksquare  &  & \\
 \vdots & & & & & & & &\blacksquare  & \\
 l_{t}  & & & & & & & &&\blacksquare  \\ \hline
\end{array}
\right\vert
    \]
  the following consistency condition  
	\begin{equation}\label{Basic_consistency}
	e_{pk}+e_{kq}+1=e_{pq}
	\end{equation}
	must hold, for all $p,k,q=1,\ldots,t$ such that $p<k<q$, so that  table $C$ is consistent if all the $\frac{t(t-1)(t-2)}{6}$ equalities are satisfied. 
\end{condition}

\noindent Let us provide the justification of consistent Condition \ref{Precies Consistency Condition}. As already observed in the section above and formally described in \cite{CorrenteGrecoSlowinski2016}, putting together eqs. (\ref{AlphaValue}) and (\ref{BlankCards}) and fixed the levels $l_{p}$ and $l_{k}$, with $p<k$, we have that 

\begin{equation}\label{FirstEqConsistency}
u(l_{k})=u(l_p)+\left(\displaystyle\sum_{r=p}^{k-1}\left(e_{r}+1\right)\right)\alpha.
\end{equation}

\noindent Supposing that among the (non-consecutive) levels $l_{p}$ and $l_{k}$ there are $e_{pk}$ blank cards, it should be true that 

\begin{equation}\label{SecondEqConsistency}
u(l_{k})=u(l_p)+(e_{pk}+1)\alpha.
\end{equation}

\noindent From eqs. (\ref{FirstEqConsistency}) and (\ref{SecondEqConsistency}), we get

\[
u(l_p)+\left(\displaystyle\sum_{r=p}^{k-1}\left(e_{r}+1\right)\right)\alpha=u(l_p)+(e_{pk}+1)\alpha\Leftrightarrow \displaystyle\sum_{r=p}^{k-1}\left(e_{r}+1\right)=e_{pk}+1
\]

\noindent and, equivalently, 

\begin{equation}\label{ThirdEqConsistency}
e_{pk}=\displaystyle\sum_{r=p}^{k-1}\left(e_{r}+1\right)-1.
\end{equation}

\noindent Applying eq. (\ref{ThirdEqConsistency}) and considering the levels $l_{p}, l_{k}, l_{q}$ such that $p<k<q,$ we have 

\begin{eqnarray}
e_{kq}&=&\displaystyle\sum_{r=k}^{q-1}\left(e_{r}+1\right)-1, \label{Eq2} \\
e_{pq}&=&\displaystyle\sum_{r=p}^{q-1}\left(e_{r}+1\right)-1. \label{Eq3} 
\end{eqnarray}

It is therefore easy to observe that 

\begin{eqnarray*}
e_{pq}&=&\displaystyle\sum_{r=p}^{q-1}\left(e_{r}+1\right)-1=\sum_{r=p}^{k-1}\left(e_{r}+1\right)+\sum_{r=k}^{q-1}\left(e_{r}+1\right)-1=\\
&=&\displaystyle \left(\sum_{r=p}^{k-1}\left(e_{r}+1\right)-1\right)+\left(\sum_{r=k}^{q-1}\left(e_{r}+1\right)-1\right)+2-1=e_{pk}+e_{kq}+1
\end{eqnarray*}

\noindent where the last equality is a consequence of eqs. (\ref{ThirdEqConsistency}) and (\ref{Eq2}). 

Let us observe that the application of the classical DCM approach involves the knowledge of the blank cards to be inserted between two successive levels, that are the values $e_{pq}$ with $p=1,\ldots,t-1$ and $q=p+1$ in the comparison table above. From these values, considering the consistency condition defined above, the whole comparison table can be filled. 

For example, let us suppose we have 5 different levels $l_{1},\ldots,l_{5}$ and that the DM applied the classical DCM approach specifying the values shown in the comparison table below

    \[
\left\vert
\begin{array}{c|ccccc}\hline
         & l_{1} & l_2 & l_3 & l_4 & l_5 \\ \hline
 l_{1}  & \blacksquare & \textcolor{red}{1} &  &  & \\
 l_{2}  & & \blacksquare & \textcolor{red}{0} & &  \\
 l_{3}  & & & \blacksquare & \textcolor{red}{3} &  \\
 l_{4}  & & & & \blacksquare & \textcolor{red}{2}  \\
 l_{5}  & & & & & \blacksquare  \\
 \hline
\end{array}
\right\vert
    \]

\noindent where $1$ is the number of blank cards that should be included between $l_1$ and $l_2$, while $2$ is the number of blank cards that should be included between $l_4$ and $l_5$. Consequently, taking into account the consistency Condition \ref{Precies Consistency Condition} defined above, we can fill all the other values in the comparison table as follows.

    \[
\left\vert
\begin{array}{c|ccccc}\hline
         & l_{1} & l_2 & l_3 & l_4 & l_5 \\ \hline
 l_{1}  & \blacksquare & 1 & \textcolor{red}{2} & \textcolor{red}{6} & \textcolor{red}{9}\\
 l_{2}  & & \blacksquare & 0 & \textcolor{red}{4} & \textcolor{red}{7} \\
 l_{3}  & & & \blacksquare & 3 & \textcolor{red}{6} \\
 l_{4}  & & & & \blacksquare & 2  \\
 l_{5}  & & & & & \blacksquare  \\
 \hline
\end{array}
\right\vert
    \]

\noindent For example, $e_{13}=e_{12}+e_{23}+1=1+0+1=2$, while $e_{25}=e_{23}+e_{35}+1=0+6+1=7$. 

\subsection{Different ways of representing the preference information}
\begin{figure}[htpb!]
        \centering
        \includegraphics[width=14cm, height=5cm]{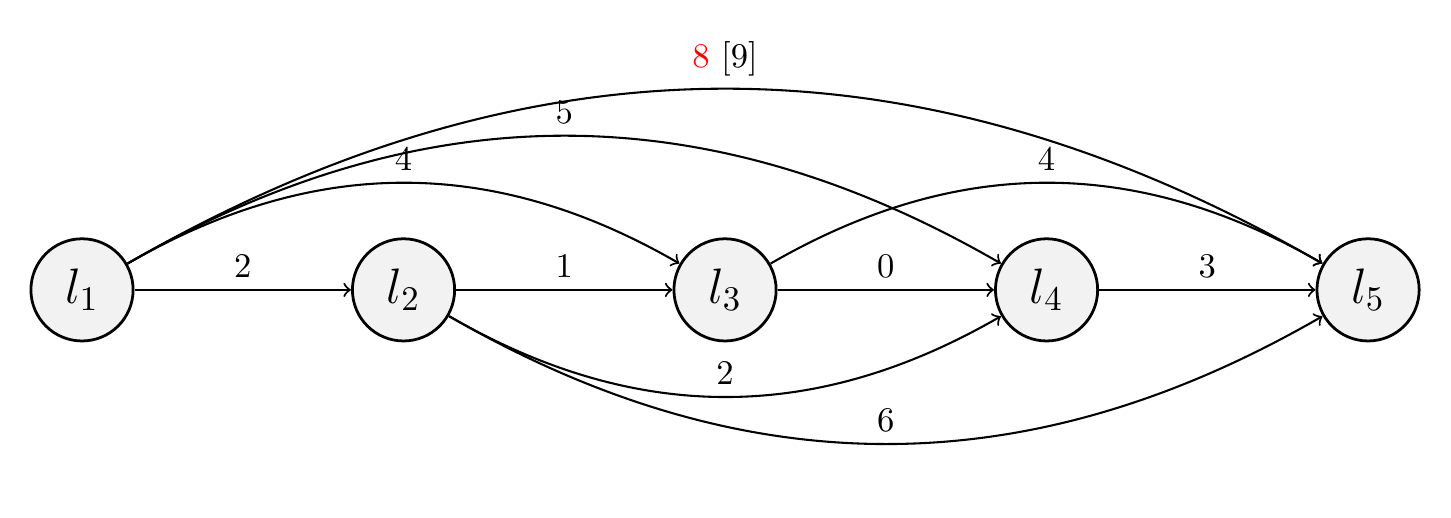}
        \caption{Graph of transitions}
        \label{Graph of transitions}
    \end{figure}

\begin{figure}[htpb!]
        \centering
        \includegraphics[width=8cm, height=10cm]{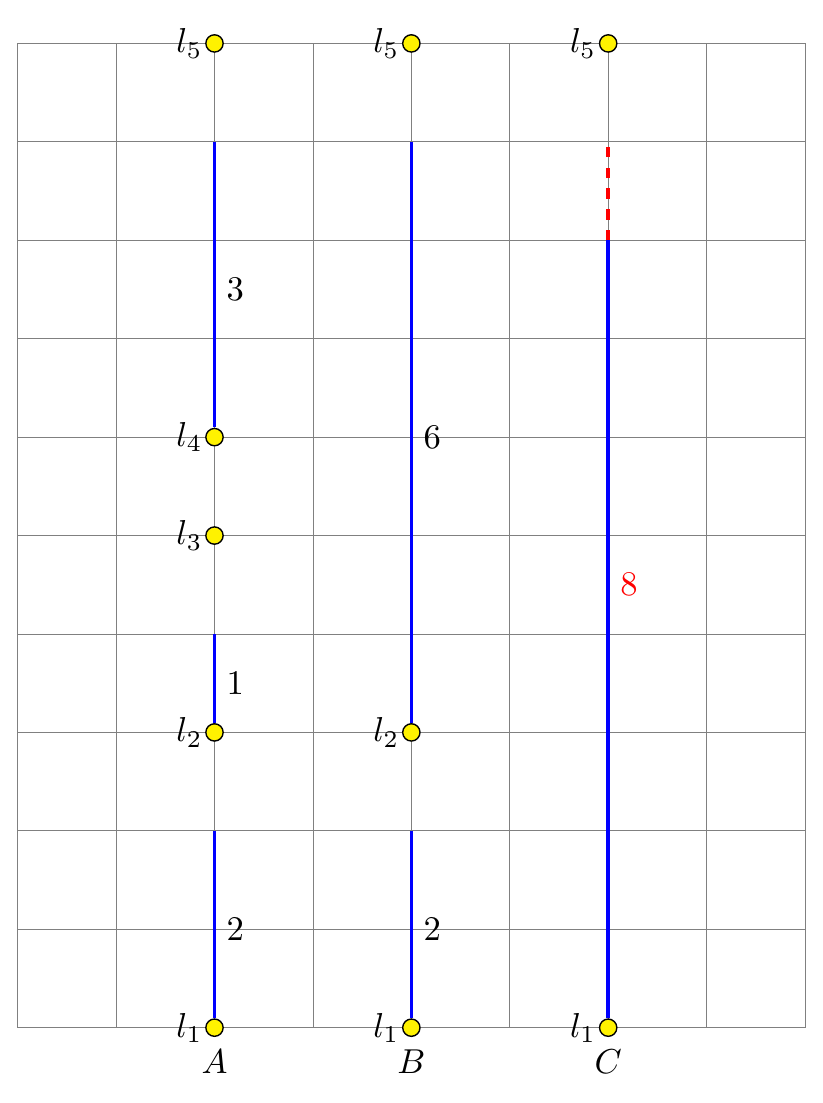}
        \caption{Bar graphs or lines}
        \label{Graph_2 of transitions}
    \end{figure}

\noindent The preference and value information collected in a DCM based comparison table can be represented in different forms as for example valued directed graphs and line bars. More in detail, given a comparison table $C$, such as for example Table \ref{comparisontable}
\begin{itemize}[label={--}]
	\item in a valued directed graph \citep{Roberts1979}, a node is assigned to each level $l_p$, an arc is linking pairs of nodes $(l_p,l_q)$ with $p>q$ and the value $e_{pq}$, that is, the number of blank cards is assigned to the corresponding arc as shown in Figure \ref{Graph of transitions},
	\item using line bars, $e_{pq}$ difference of appreciations between levels $l_p$ and $l_q$ are represented as segments having a length with a measure $e_{pq}+1$ as shown in Figure \ref{Graph_2 of transitions}.   
\end{itemize}

\begin{table}[ht]
\caption{Deck of cards comparison table} % title of Table
\centering % used for centering table
\begin{tabular}{c|ccccc} % centered columns (4 columns)
\hline %inserts double horizontal lines
& $l_{1}$ & $l_{2}$ & $l_{3}$ & $l_{4}$ & $l_{5}$ \\ \hline
$l_{1}$ & $\blacksquare$ & 2 & 4 & 5 & \textcolor{red}{8} \\
$l_{2}$ & & $\blacksquare$ & 1 & 2 & 6 \\
$l_{3}$ & &   & $\blacksquare$ & 0 & 4 \\
$l_{4}$ & &   &   & $\blacksquare$ & 3 \\
$l_{5}$ & &   &   &   & $\blacksquare$ \\ \hline
\end{tabular}
\label{comparisontable} % is used to refer this table in the text
\end{table}
--------------- 
%the  information to the DM and the interaction protocols are of utmost importance for gathering the richest pieces of preference information. There are several tools we can use to interact with the DM: tables, graphs, line bars (these are presented below, and they are quite easy to interpret). 
%\begin{table}[ht]
%\caption{Deck of card based comparison table}
%\[
%\left\vert
%\begin{array}{c|ccccc}\hline
       %& l_{1} & l_{2} & l_{3} & l_{4} & l_{5} \\ \hline
%l_{1} & \blacksquare & 2 & 4 & 5 & \textcolor{red}{8} \\
%l_{2} & & \blacksquare & 1 & 2 & 6 \\
%l_{3} & &   & \blacksquare & 0 & 4 \\
%l_{4} & &   &   & \blacksquare & 3 \\
%l_{5} & &   &   &   & \blacksquare \\ \hline
%\end{array}
%\right\vert
%\]
%\end{tabular}
%\label{table:nonlin} 
%\end{table}

%\noindent Graph (Figure \ref{Graph of transitions}) represents the same preference information as the table before. 

\vfill\newpage 

%\noindent Analogously, bar graphs or lines (Figure \ref{Graph_2 of transitions}) are tools that can be used to represent the preference information and render it easy to he DM.  

%%%%%%%%%%%%%%%%%%%%%%%%%%%%%%%%%%%%%%%%%%%%%%%%%%%%%%%%%%%%%%%%%%%%%%%%%%%%%%%%%%%%%%%%%%%%%%
\section{Detecting and correcting inconsistency in pairwise comparison tables based on the deck of cards method }\label{RestoringInconsistency}%%
%%%%%%%%%%%%%%%%%%%%%%%%%%%%%%%%%%%%%%%%%%%%%%%%%%%%%%%%%%%%%%%%%%%%%%%%%%%%%%%%%%%%%%%%%%%%%%

\noindent In the behavioural aspects of decision-making it is important to develop tools for dealing with inconsistent judgments and ways of interacting with the DM to overcome, when possible, such inconsistency. In case of inconsistent judgments, the comparison table provided by the DM does not satisfy the consistency Condition \ref{Precies Consistency Condition}. Therefore, in this section we present some tools to detect and correct the inconsitency. These are mainly mixed-integer linear programming (MILP) based methods, which are  used to check for the consistency of the preference and value information provided by the DM. If this information is not consistent, the MILP program will suggest the minimum number of modifications necessary to restore the consistency. An example will be provided in this sections as well as some final considerations about the decision-making and decision-aiding behavioural aspects of our DCM approach. 

To check if the pairwise comparison table $C$ is consistent, that is Condition \ref{Precies Consistency Condition} is satisfied, we shall proceed in the following way. For each ordered pair of levels $(l_p,l_q)$ such that $p<q$, we define three different variables $\delta^{+}_{pq}$, $\delta^{-}_{pq}$ and $\overline{e}_{pq}$ such that: 

\[
\overline{e}_{pq}=e_{pq}+\delta^{+}_{pq}-\delta^{-}_{pq}.
\]

Their meaning is as follows: 

\begin{itemize}[label={--}]
\item $\delta^{-}_{pq}$ is a non-negative integer number, representing how many blank cards should be subtracted from $e_{pq}$, that is, how much the difference in appreciation $e_{pq}$ has to be reduced to make the judgments consistent;
\item $\delta^{+}_{pq}$ is a non-negative integer number, representing how many blank cards should be added to $e_{pq}$, that is, how much the difference in appreciation $e_{pq}$ has to be increased to make the judgments consistent;
\item $\overline{e}_{pq}$ is the \textit{new} number of blank cards that should be included between levels $l_p$ and $l_q$, that is the difference in appreciation between $l_p$ and $l_q$ to make consistent the comparison table. 
\end{itemize}

Moreover, for each pair of levels $(l_p,l_q)$, a binary variable $y_{pq}$ is defined to check if the starting evaluation $e_{pq}$ has to be modified or not. Let $P$ denote the set of all feasible ordered pairs $(p,q)$ in the comparison table, that is $P=\{(p,q): p,q=1,\ldots,t \;\mbox{and}\; p<q\}$.

To check if the comparison table provided by the DM is consistent, one has to solve the following MILP problem: 

\begin{subequations}
\begin{gather}
\hspace{-5cm} y^{*}=\argmin z(y)=\displaystyle\sum_{(p,q) \; \in \; P}y_{pq}  \nonumber \\
\intertext{\hspace{3.75cm}subject to:}
 \hspace{2.5cm}  \overline{e}_{pk}+\overline{e}_{kq}+1=\overline{e}_{pq},  \qquad  (p,k),\, (k,q), \; (p,q) \in P, \label{Consistence Constraint_1} \\ 
 \hspace{0.5cm}   e_{pq}+\delta_{pq}^{+}-\delta_{pq}^{-}=\overline{e}_{pq},  \qquad  (p,q)  \in P,  \label{Consistence Constraint_2} \\
 \hspace{0.5cm}   \delta^{+}_{pq}+\delta^{-}_{pq}\leqslant  My_{pq},\;\;\;\; \; \qquad (p,q)  \in P,  \label{Consistence Constraint_3} \\
  \hspace{1.5cm}  \delta_{pq}^{-},\, \delta_{pq}^{+} \in \Bbb{N}_0, \;\;\; \qquad (p,q)  \in P, \label{Consistence Constraint_4} \\
   \hspace{1.5cm}  y_{pq}\in\{0,1\},\;\;\;\; \qquad (p,q)  \in P. \label{Consistence Constraint_5} \\
\intertext{} \nonumber
\end{gather}
\end{subequations}

\vspace{-1.5cm}

In this model $M$ is a big positive number. The objective function $z(y)$, that has to be minimized, counts the number of modifications necessary to make consistent the comparison table $C$. Constraint (\ref{Consistence Constraint_1}) is the consistency condition of the new comparison table that should be fulfilled by all triples of levels $(p,k,q)$ such that $p<k<q$. Constraint (\ref{Consistence Constraint_2}) is used to link the evaluations in the starting comparison table $(e_{pq})$ to the new ones $(\overline{e}_{pq})$\footnote{Note that Constraints (\ref{Consistence Constraint_1}) and (\ref{Consistence Constraint_2}) can be converted into a single one, that is, $\overline{e}_{pk}+\overline{e}_{kq}+1=e_{pq}+\delta_{pq}^{+}-\delta_{pq}^{-}$, but, in this way, the model looses readability.}. Constraint (\ref{Consistence Constraint_3}) is used to check if the comparison $e_{pq}$ has to be modified or not. If $y_{pq}=0$, then $\delta^{+}_{pq}=\delta^{-}_{pq}=0$ and, consequently, the evaluation $e_{pq}$ has not to be modified, while if $y_{pq}=1$, then Constraint (\ref{Consistence Constraint_3}) is always satisfied and, therefore, $e_{pq}$ has to be modified. Constraints (\ref{Consistence Constraint_4}) and (\ref{Consistence Constraint_5}) are used to express the nature of the used variables. 

If $z^{*}=0$, with $z^*$ being the optimal value of the objective function $z$, then the comparison table is consistent and, therefore, no evaluation $e_{pq}$ needs to be modified. In the opposite case, the comparison table is not consistent and the evaluations for which $y^{*}_{pq}=1$ need to be modified adding $\delta_{pq}^{+}$ units to $e_{pq}$ or reducing of $\delta_{pq}^{-}$ units the same evaluation. Let us observe that $z^{*}$, the optimal value of our objective function, is the number of 1s in the vector $y^{*}$. 

If the comparison table is not consistent, i.e., $z^{*}>0$, it is reasonable to find all the possible sets of modifications, which are necessary to restore its consistency. To check for another possible solution, we need to solve exactly the same problem above with the addition of the following constraint: 

\begin{equation}
\displaystyle\sum_{\{ (p,q)\; \in P:\; y^{\ast}_{pq} = 1\}}y_{pq}\leqslant z^{\ast}-1. 
\end{equation}

 This constraint is used to forbid to obtain the sets of possible modifications found in the previous iterations. The procedure continues until the MILP problem becomes infeasible and, therefore, does not exist any other modification of the evaluations making the comparison table consistent. 
 
 In the following subsection we present an example to illustrate these concepts. 

%%%%%%%%%%%%%%%%%%%%%%%%%%%%%%%%%%%%%%%
\subsection{An illustrative example}%%%
%%%%%%%%%%%%%%%%%%%%%%%%%%%%%%%%%%%%%%%
\noindent Let us consider five different levels, $l_{1},l_{2},l_{3},l_{4},l_{5}$ ordered from the worst to the best one. Let us assume that the DM is able to provide the comparison table below. 
 
\[
\left\vert
\begin{array}{c|ccccc}\hline
         & l_{1} & l_{2} & l_{3} & l_{4} & l_{5} \\ \hline
 l_{1} & \blacksquare & 2 & 4 & 5 & 9 \\
 l_{2} & & \blacksquare  & 1 & 2 & 6 \\
 l_{3} & &   &  \blacksquare & 0 & 4 \\
 l_{4} & &   &   &  \blacksquare & 3 \\
 l_{5} & &   &   &   &  \blacksquare \\ \hline
\end{array}
\right\vert
\]

It is easy to observe that the matrix is perfectly consistent since $e_{pq}=e_{pk}+e_{kq}+1$ for all $(p,k),(k,q),(p,q)\in P$.  

Now, let us consider another example, in which the DM is able to provide the information in the following comparison table.  

\[
\left\vert
\begin{array}{c|ccccc}\hline
       & l_{1} & l_{2} & l_{3} & l_{4} & l_{5} \\ \hline
l_{1} & \blacksquare & 2 & 4 & 5 & \textcolor{red}{8} \\
l_{2} & & \blacksquare & 1 & 2 & 6 \\
l_{3} & &   & \blacksquare & 0 & 4 \\
l_{4} & &   &   & \blacksquare & 3 \\
l_{5} & &   &   &   & \blacksquare \\ \hline
\end{array}
\right\vert
\]

Solving the first MILP problem described in section \ref{RestoringInconsistency}, we obtain $z^{*}>0$ and, consequently, the comparison table provided by the DM is not consistent. In this case, we can observe that the inconsistency is present since, for example, $e_{14}+e_{45}+1=5+3+1=9\neq e_{15}=8$. This means that the intensity of the preference of $l_{5}$ over $l_{1}$ $(e_{15})$ is lower than the sum of the two intensities of preference considered in this case, that is, that one between $l_{4}$ and $l_{1}$ ($e_{14}$) and that one between $l_{5}$ and $l_{4}$ ($e_{45}$) plus one. All the possible modifications making consistent the starting comparison table together with the new consistent comparison table and the constraints to be added to obtain the new sets of modifications are listed in the following: 
 
\begin{enumerate}
\item Restoring consistency requires only one change: $e_{15}\leftarrow 9$ (instead of 8). Since $y_{15}^\ast =  1$, then the following constraint must be added to the MILP problem: $y_{15} \leqslant 0$.  

\[
\left\vert
\begin{array}{c|ccccc}\hline
         & l_{1} & l_{2} & l_{3} & l_{4} & l_{5} \\ \hline
 l_{1} & \blacksquare & 2 & 4 & 5 & \textcolor{red}{9} \\
 l_{2} & & \blacksquare  & 1 & 2 & 6 \\
 l_{3} & &   & \blacksquare  & 0 & 4 \\
 l_{4} & &   &   & \blacksquare  & 3 \\
 l_{5} & &   &   &   &  \blacksquare \\ \hline
\end{array}
\right\vert
\]

\item Restoring consistency requires three changes: $e_{25}\leftarrow 5$ (instead of 6); $e_{35}\leftarrow 3$ (instead of 4); and,  $e_{45}\leftarrow 2$ (instead of 3). Since $y_{25}^\ast = y_{35}^\ast = y_{45}^\ast = 1$, then the following constraint must be added to the MILP problem: $y_{25} + y_{35} + y_{45} \leqslant 2$.

\[
\left\vert
\begin{array}{c|ccccc}\hline
       & l_{1} & l_{2} & l_{3} & l_{4} & l_{5} \\ \hline
l_{1} & \blacksquare & 2 & 4 & 5 & 8 \\
l_{2} & & \blacksquare  & 1 & 2 & \textcolor{red}{5} \\
l_{3} & &   &  \blacksquare & 0 & \textcolor{red}{3} \\
l_{4} & &   &   & \blacksquare  & \textcolor{red}{2} \\
l_{5} & &   &   &   &  \blacksquare \\ \hline 
\end{array}
\right\vert
\] 

 \item Restoring consistency requires also three changes: $e_{12}\leftarrow 1$ (instead of 2); $e_{13}\leftarrow 3$ (instead of 4); and,   $e_{14}\leftarrow 4$ (instead of 5). Since $y_{12}^\ast = y_{13}^\ast = y_{14}^\ast = 1$, then the following constraint must be added to the MILP problem: $y_{12} + y_{13} + y_{14} \leqslant 2$.

\[ 
\left\vert
\begin{array}{c|ccccc}\hline
       & l_{1} & l_{2} & l_{3} & l_{4} & l_{5} \\ \hline
l_{1} & \blacksquare & \textcolor{red}{1} & \textcolor{red}{3} & \textcolor{red}{4} & 8 \\
l_{2} & &  \blacksquare & 1 & 2 & 6 \\
l_{3} & &   &  \blacksquare & 0 & 4 \\
l_{4} & &   &   & \blacksquare  & 3 \\
l_{5} & &   &   &   &  \blacksquare \\ \hline
\end{array}
\right\vert
\]

\item Restoring consistency requires now five changes: $e_{13}\leftarrow 3$ (instead of 4); $e_{14}\leftarrow 4$ (instead of 5); $e_{23}\leftarrow 0$ (instead of 1); $e_{24}\leftarrow 1$ (instead of 2); and,   $e_{25}\leftarrow 5$ (instead of 6). Since $y_{13}^\ast = y_{14}^\ast = y_{23}^\ast = y_{24}^\ast = y_{25}^\ast = 1$, then the following constraint must be added to the MILP problem: $y_{13} + y_{14} + y_{23} + y_{24} + y_{25} \leqslant 4$.

\[
\left\vert
\begin{array}{c|ccccc}\hline
       & l_{1} & l_{2} & l_{3} & l_{4} & l_{5} \\ \hline
l_{1} & \blacksquare & 2 & \textcolor{red}{3} & \textcolor{red}{4} & 8 \\
l_{2} & & \blacksquare  & \textcolor{red}{0} & \textcolor{red}{1} & \textcolor{red}{5} \\
l_{3} & &   &  \blacksquare & 0 & 4 \\
l_{4} & &   &   & \blacksquare  & 3 \\
l_{5} & &   &   &   &  \blacksquare \\ \hline
\end{array}
\right\vert
\]

\item Restoring consistency requires also five changes: $e_{13}\leftarrow 3$ (instead of 4); $e_{23}\leftarrow 0$ (instead of 1); $e_{25}\leftarrow 5$ (instead of 6); $e_{34}\leftarrow 1$ (instead of 0); and,   $e_{45}\leftarrow 2$ (instead of 3). Since $y_{13}^\ast = y_{23}^\ast = y_{25}^\ast = y_{34}^\ast = y_{45}^\ast = 1$, then the following constraint must be added to the MILP problem: $y_{13} + y_{23} + y_{25} + y_{34} + y_{45} \leqslant 4$.

\[
\left\vert
\begin{array}{c|ccccc}\hline
       & l_{1} & l_{2} & l_{3} & l_{4} & l_{5} \\ \hline
l_{1} & \blacksquare & 2 & \textcolor{red}{3} & 5 & 8 \\
l_{2} & &  \blacksquare & \textcolor{red}{0} & 2 & \textcolor{red}{5} \\
l_{3} & &   & \blacksquare  & \textcolor{red}{1} & 4 \\
l_{4} & &   &   & \blacksquare  & \textcolor{red}{2} \\
l_{5} & &   &   &   &  \blacksquare \\ \hline
\end{array}
\right\vert
\] 
  
\item Restoring consistency requires again five changes: $e_{12}\leftarrow 1$ (instead of 2); $e_{23}\leftarrow 2$ (instead of 1); $e_{24}\leftarrow 3$ (instead of 2); $e_{35}\leftarrow 3$ (instead of 4); and,   $e_{45}\leftarrow 2$ (instead of 3). Since $y_{12}^\ast = y_{23}^\ast = y_{24}^\ast = y_{35}^\ast = y_{45}^\ast = 1$, then the following constraint must be added to the MILP problem: $y_{12} + y_{23} + y_{24} + y_{35} + y_{45} \leqslant 4$.

\[
\left\vert
\begin{array}{c|ccccc}\hline
       & l_{1} & l_{2} & l_{3} & l_{4} & l_{5} \\ \hline 
l_{1} & \blacksquare & \textcolor{red}{1} & 4 & 5 & 8 \\
l_{2} & & \blacksquare  & \textcolor{red}{2} & \textcolor{red}{3} & 6 \\
l_{3} & &   &  \blacksquare & 0 & \textcolor{red}{3} \\
l_{4} & &   &   &  \blacksquare & \textcolor{red}{2} \\
l_{5} & &   &   &   &  \blacksquare \\ \hline
\end{array}
\right\vert
\]

\item Restoring consistency requires once more five changes: $e_{12}\leftarrow 1$ (instead of 2); $e_{13}\leftarrow 3$ (instead of 4); $e_{24}\leftarrow 3$ (instead of 2); $e_{34}\leftarrow 1$ (instead of 0); and,   $e_{45}\leftarrow 2$ (instead of 3). Since $y_{12}^\ast = y_{13}^\ast = y_{24}^\ast = y_{34}^\ast = y_{45}^\ast = 1$, then the following constraint must be added to the MILP problem: $y_{12} + y_{13} + y_{24} + y_{34} + y_{45} \leqslant 4$.

\[
\left\vert
\begin{array}{c|ccccc}\hline 
       & l_{1} & l_{2} & l_{3} & l_{4} & l_{5} \\ \hline
l_{1} & \blacksquare & \textcolor{red}{1} & \textcolor{red}{3} & 5 & 8 \\
l_{2} & & \blacksquare  & 1 & \textcolor{red}{3} & 6 \\
l_{3} & &   &  \blacksquare  & \textcolor{red}{1} & 4 \\
l_{4} & &   &   & \blacksquare  & \textcolor{red}{2} \\
l_{5} & &   &   &   &  \blacksquare \\ \hline
\end{array}
\right\vert
\]

\end{enumerate}

After that, the MILP problem becomes infeasible and, therefore, there is not any other possibility to restore the consistency of the starting comparison table.

%%%%%%%%%%%%%%%%%%%%%%%%%%%%%%%%%%%%%%%%%%%%%%%%%%%%%%%%%%%%%%
\section{Some further developments}\label{sec:ImprovedDeck}%%%
%%%%%%%%%%%%%%%%%%%%%%%%%%%%%%%%%%%%%%%%%%%%%%%%%%%%%%%%%%%%%%
\noindent The section is devoted to some further developments of the method presented in the previous sections. It comprises the cases of missing data and the case in which information is provided in a imprecise way through intervals. We also present some approaches to deal with missing and imprecise preference judgments. 

%%%%%%%%%%%%%%%%%%%%%%%%%%%%%%%%%%%%%%%%%%%%%%%%%
\subsection{Missing values}\label{sec:Missing}%%%
%%%%%%%%%%%%%%%%%%%%%%%%%%%%%%%%%%%%%%%%%%%%%%%%%
\noindent In this section, we consider the case in which some of the comparisons in the table are missing. This means that the DM was not able to provide any value for the number of cards to be inserted between two different levels. Using the notation introduced in section \ref{RestoringConsistency}, we will describe a procedure that will be composed of the following main steps: 

\begin{enumerate}
\item check if there exists at least one set of values to be assigned to the missing comparisons so that the consistency Condition \ref{Precies Consistency Condition} is satisfied. If this is not the case, the inconsistent judgments will be highlighted to be shown to the DM and, eventually, revised,
\item in case the information provided by the DM is consistent, the procedure will show how to obtain different sets of values compatible with such information.
\end{enumerate} 

To check if there exists at least one set of values to be assigned to the missing judgments so that the preference information provided by the DM is consistent, we have to solve the following MILP problem denoted in the following by $MILP-{M_1}$,   

\begin{subequations}
\begin{gather}
\hspace{-5cm} y^{*}=\argmin z(y)=\displaystyle\sum_{(p,q) \; \in \; P}y_{pq}  \nonumber \\
\intertext{\hspace{3.75cm}subject to:}
 \hspace{2.5cm}  \overline{e}_{pk}+\overline{e}_{kq}+1=\overline{e}_{pq},  \qquad  (p,k), \, (k,q),\, (p,q) \in P, \label{eq: Constraint_1} \\ 
 \hspace{0.5cm}   e_{pq}+\delta_{pq}^{+}-\delta_{pq}^{-}=\overline{e}_{pq},  \qquad  (p,q)  \in P^{DM},  \label{eq: Constraint_2} \\
 \hspace{0.5cm}   \delta^{+}_{pq}+\delta^{-}_{pq}\leqslant  My_{pq},\;\;\;\; \; \qquad (p,q)  \in P^{DM},  \label{eq: Constraint_3} \\
  \hspace{1.5cm}  \delta_{pq}^{-},\, \delta_{pq}^{+} \in \Bbb{N}_0, \;\;\; \qquad (p,q)  \in P^{DM}, \label{eq: Constraint_4} \\
   \hspace{1.5cm}  y_{pq}\in\{0,1\},\;\;\;\; \qquad (p,q)  \in P^{DM} \label{eq: Constraint_5} \\
   \hspace{1.5cm}  \overline{e}_{pq}\in\Bbb{N}_{0},\;\;\;\; \qquad (p,q)  \in P \label{eq: Constraint_6} \\
\intertext{} \nonumber
\end{gather}
\end{subequations}

\vspace{-1.5cm}

\noindent where $P^{DM}\subseteq P$ is the set of pairs of levels for which the DM was able to provide some preference information. Let us underline that, differently from the MILP problem defined in the previous section, the auxiliary variables $\delta_{pq}^{+},\delta_{pq}^{-}$ and $y_{pq}$ are now defined for the pairs $(p,q)\in P^{DM}$ only, i.e., the pairs of levels for which the DM expressed a preference information. 

If $z^{*}>0$, then the comparisons provided by the DM are not consistent and, therefore, need to be revised. In particular, the pairs $(p,q)\in P^{DM}$ for which $y_{pq}^{*}=1$ are those needing to be modified to restore the consistency. 

If $z^{*}=0$, then the few comparisons provided by the DM are consistent and, consequently, $e_{pq}=\overline{e}^{*}_{pq}$ for all $(p,q)\in P^{DM}$, while the values $\overline{e}^{*}_{pq}$ with $(p,q)\in P\setminus P^{DM}$ will be the values to be assigned to the missing information. In this case, one has to check if the found solution is unique or if there exists some other solution. \\
Denoting by $\overline{e}^{*}_{pq}$ for all $(p,q)\in P\setminus P^{DM}$ the values assigned to the missing information by the solution of the previous problem, to check for another possible solution, one has to solve the MILP problem, denoted by $MILP-{M_2}$, obtained adding to $MILP-{M_1}$ the following constraints:

\begin{subequations}
\begin{gather}
 \hspace{0.5cm}  \overline{e}_{pq}\leqslant \left(\overline{e}^{*}_{pq}-1\right)+My^{1}_{pq},  \qquad  (p,q)\in P\setminus P^{DM}, \label{eq: Constraint_1_Further} \\ 
 \hspace{0.5cm}  \overline{e}_{pq}\geqslant \left(\overline{e}^{*}_{pq}+1\right)-My^{2}_{pq},  \qquad  (p,q)\in P\setminus P^{DM}, \label{eq: Constraint_2_Further} \\ 
 \hspace{0.5cm}   \displaystyle\sum_{(p,q)\in P\setminus P^{DM}}\left(y^{1}_{pq}+y^{2}_{pq}\right)\leqslant 2\vert P\setminus P^{DM}\vert-1,\;\;\;\; \;   \label{eq: Constraint_3_Further} \\
   \hspace{1.5cm}  y^{1}_{pq},y^{2}_{pq}\in\{0,1\},\;\;\;\; \qquad (p,q)  \in P\setminus P^{DM}. \label{eq: Constraint_4_Further} 
\end{gather}
\end{subequations}

 If $y_{pq}^{1}=1$ or $y_{pq}^{2}=1$, the corresponding constraints $(\ref{eq: Constraint_1_Further})$ and $(\ref{eq: Constraint_2_Further})$ are always satisfied and, consequently, $\overline{e}_{pq}=\overline{e}^{*}_{pq}$. Constraint $(\ref{eq: Constraint_3_Further})$ ensures, therefore, that holds at least one between $\overline{e}_{pq}\geqslant \overline{e}^{*}_{pq}+1$ and $\overline{e}_{pq}\leqslant \overline{e}^{*}_{pq}-1$ for at least one pair $(p,q)\in P\setminus P^{DM}$.  

If the set of constraints in $MILP-M_{2}$ is feasible and $z^{*}=0,$ then we were able to obtain another solution, while, in the opposite case, the found solution is the unique one. Proceeding in an iterative way one can obtain other possible solutions. 

%%%%%%%%%%%%%%%%%%%%%%%%%%%%%%%%%%%%%%%%%%%%%%%%%%%%%%%%%%%%%%%%%%%%%%
\subsubsection{An illustrative example}\label{MissingValuesExample}%%%
%%%%%%%%%%%%%%%%%%%%%%%%%%%%%%%%%%%%%%%%%%%%%%%%%%%%%%%%%%%%%%%%%%%%%%
\noindent Let us assume that the DM was able to provide the preference information shown in the following table where the positions containing ``$?$" denote some missing values since the DM was not able or was not sure about that comparison.  

\[
\left\vert
\begin{array}{c|ccccc}\hline
           & l_{1} & l_{2} & l_{3} & l_{4} & l_{5} \\ \hline
    l_{1}  &\blacksquare & {?} & {?} & {?} & {?} \\
    l_{2}  & &\blacksquare & {?} & {?}  & 6   \\
    l_{3}  & & &\blacksquare & 0  & 5      \\
    l_{4}  & & & &\blacksquare & 3       \\
    l_{5}  & & & & &\blacksquare         \\  \hline
\end{array}
\right\vert
\]

Solving $MILP-{M_1}$, we obtain $z^{*}=1$ and, in particular, $y^{*}_{35}=1$ meaning that the value corresponding to the comparison between $l_3$ and $l_5$ is not consistent with the other information provided by the DM. Indeed, considering the consistency Condition \ref{Precies Consistency Condition}, we should have $e_{34}+e_{45}+1=e_{35}$ being not fulfilled since $0+3+1\neq 5$. Solving the problem, we get $\delta^{-}_{35}=1$ meaning that if the evaluation $e_{35}=5\rightarrow 4$, then the consistency is again restored. Doing the suggested replacement we get the new comparison table 

\[
\left\vert
\begin{array}{c|ccccc}\hline
           & l_{1} & l_{2} & l_{3} & l_{4} & l_{5} \\ \hline
    l_{1}  &\blacksquare & \textcolor{red}{2} & \textcolor{red}{4} & \textcolor{red}{5} & \textcolor{red}{9} \\
    l_{2}  & &\blacksquare & \textcolor{red}{1} & \textcolor{red}{2}  & 6   \\
    l_{3}  & & &\blacksquare & 0  & 4      \\
    l_{4}  & & & &\blacksquare & 3       \\
    l_{5}  & & & & &\blacksquare         \\  \hline
\end{array}
\right\vert
\]

At this point, to check if there exists another possible solution, we avoid to obtain again the same solution adding the following constraints to $MILP-{M_1}$: 

\begin{subequations}
\begin{gather}
 \hspace{0.5cm}  \overline{e}_{12}\leqslant \left(2-1\right)+My^{1}_{12}, \\ 
 \hspace{0.5cm}  \overline{e}_{12}\geqslant \left(2+1\right)-My^{2}_{12}, \\ 
 \hspace{0.5cm}  \overline{e}_{13}\leqslant \left(4-1\right)+My^{1}_{13}, \\ 
 \hspace{0.5cm}  \overline{e}_{13}\geqslant \left(4+1\right)-My^{2}_{13}, \\
 \hspace{0.5cm}  \overline{e}_{14}\leqslant \left(5-1\right)+My^{1}_{14}, \\ 
 \hspace{0.5cm}  \overline{e}_{14}\geqslant \left(5+1\right)-My^{2}_{14}, \\ 
 \hspace{0.5cm}  \overline{e}_{15}\leqslant \left(9-1\right)+My^{1}_{15}, \\ 
 \hspace{0.5cm}  \overline{e}_{15}\geqslant \left(9+1\right)-My^{2}_{15}, \\ 
 \hspace{0.5cm}  \overline{e}_{23}\leqslant \left(1-1\right)+My^{1}_{23}, \\ 
 \hspace{0.5cm}  \overline{e}_{23}\geqslant \left(1+1\right)-My^{2}_{23}, \\
 \hspace{0.5cm}  \overline{e}_{24}\leqslant \left(2-1\right)+My^{1}_{24}, \\ 
 \hspace{0.5cm}  \overline{e}_{24}\geqslant \left(2+1\right)-My^{2}_{24}, \\ 
 \hspace{0.5cm}   \displaystyle\sum_{(p,q)\in \left\{(1,2),(1,3),(1,4),(1,5),(2,3),(2,4)\right\}}\left(y^{1}_{pq}+y^{2}_{pq}\right)\leqslant 11,\;\;\;\; \; \\
   \hspace{1.5cm}  y^{1}_{pq},y^{2}_{pq}\in\{0,1\},\;\;\;\; \qquad (p,q)  \in \left\{(1,2),(1,3),(1,4),(1,5),(2,3),(2,4)\right\}. 
\end{gather}
\end{subequations}

Solving the new MILP problem ($MILP-{M_2}$), we get $z^{*}=0$ and, therefore, there exists another possible solution being shown in the following comparison table

\[
\left\vert
\begin{array}{c|ccccc}\hline
           & l_{1} & l_{2} & l_{3} & l_{4} & l_{5} \\ \hline
    l_{1}  &\blacksquare & \textcolor{red}{0} & \textcolor{red}{2} & \textcolor{red}{3} & \textcolor{red}{7} \\
    l_{2}  & &\blacksquare & \textcolor{red}{1} & \textcolor{red}{2}  & 6   \\
    l_{3}  & & &\blacksquare & 0  & 4      \\
    l_{4}  & & & &\blacksquare & 3       \\
    l_{5}  & & & & &\blacksquare         \\  \hline
\end{array}
\right\vert
\]

Let us observe that in the new comparison table, the values corresponding to $(2,3)$ and $(2,4)$ have not been changed while the others have been modified. One could continue by applying the same procedure to check for other possible sets of values. 

Now, let us assume that the preference information provide by the DM is contained in the comparison table below: 

\[
\left\vert
\begin{array}{c|ccccc}\hline
           & l_{1} & l_{2} & l_{3} & l_{4} & l_{5} \\ \hline
    l_{1}  &\blacksquare & {?} & 4 & {?} & 9 \\
    l_{2}  & &\blacksquare & {?} & {?}  & 6   \\
    l_{3}  & & &\blacksquare & {?}  & {?}      \\
    l_{4}  & & & &\blacksquare & 3       \\
    l_{5}  & & & & &\blacksquare         \\  \hline
\end{array}
\right\vert
\]

Solving $MILP-{M_1}$, we obtain $z^{*}=0$ meaning that there is at least one comparison table compatible with the few information given by the DM and the obtained solution is the following: 

\[
\left\vert
\begin{array}{c|ccccc}\hline
           & l_{1} & l_{2} & l_{3} & l_{4} & l_{5} \\ \hline
    l_{1}  &\blacksquare & \textcolor{red}{2} & 4 & \textcolor{red}{5} & 9 \\
    l_{2}  & &\blacksquare & \textcolor{red}{1} & \textcolor{red}{2}  & 6   \\
    l_{3}  & & &\blacksquare & \textcolor{red}{0}  & \textcolor{red}{4}      \\
    l_{4}  & & & &\blacksquare & 3       \\
    l_{5}  & & & & &\blacksquare         \\  \hline
\end{array}
\right\vert
\]

To check for the existence of another comparison table, we have to solve $MILP-{M_2}$ obtained by $MILP-{M_1}$ with the addition of the following constraints: 

\begin{subequations}
\begin{gather}
 \hspace{0.5cm}  \overline{e}_{12}\leqslant \left(2-1\right)+My^{1}_{12}, \\ 
 \hspace{0.5cm}  \overline{e}_{12}\geqslant \left(2+1\right)-My^{2}_{12}, \\ 
 \hspace{0.5cm}  \overline{e}_{14}\leqslant \left(5-1\right)+My^{1}_{14}, \\ 
 \hspace{0.5cm}  \overline{e}_{14}\geqslant \left(5+1\right)-My^{2}_{14}, \\ 
 \hspace{0.5cm}  \overline{e}_{23}\leqslant \left(1-1\right)+My^{1}_{23}, \\ 
 \hspace{0.5cm}  \overline{e}_{23}\geqslant \left(1+1\right)-My^{2}_{23}, \\
 \hspace{0.5cm}  \overline{e}_{24}\leqslant \left(2-1\right)+My^{1}_{24}, \\ 
 \hspace{0.5cm}  \overline{e}_{24}\geqslant \left(2+1\right)-My^{2}_{24}, \\ 
 \hspace{0.5cm}  \overline{e}_{34}\leqslant \left(0-1\right)+My^{1}_{34}, \\ 
 \hspace{0.5cm}  \overline{e}_{34}\geqslant \left(0+1\right)-My^{2}_{34}, \\
 \hspace{0.5cm}  \overline{e}_{35}\leqslant \left(4-1\right)+My^{1}_{35}, \\ 
 \hspace{0.5cm}  \overline{e}_{35}\geqslant \left(4+1\right)-My^{2}_{35}, \\
 \hspace{0.5cm}   \displaystyle\sum_{(p,q)\in \left\{(1,2),(1,4),(2,3),(2,4),(3,4),(3,5)\right\}}\left(y^{1}_{pq}+y^{2}_{pq}\right)\leqslant 11,\;\;\;\; \; \\
   \hspace{1.5cm}  y^{1}_{pq},y^{2}_{pq}\in\{0,1\},\;\;\;\; \qquad (p,q)  \in \left\{(1,2),(1,4),(2,3),(2,4),(3,4),(3,5)\right\}. 
\end{gather}
\end{subequations}

In this case, we obtain $z^{*}>0$ meaning that we cannot obtain another consistent comparison table without modifying the preference information provided by the DM. Therefore, the found comparison table is unique. 

\subsection{The imprecise judgments in the form of intervals}\label{sec:Imprecise}%%%%
%%%%%%%%%%%%%%%%%%%%%%%%%%%%%%%%%%%%%%%%%%%%%%%%%%%%%%%%%%%%%%%%%%%%%%%%%%%%%%%%%%%%%%
\noindent Now we shall develop the whole methodology corresponding to the case in which the DM is not able to give a precise value to the number of blank cards added between the levels $l_p$ and $l_q$ but he prefers to express an imprecise evaluation represented by an interval of possible values. For this reason, the evaluation $e_{pq}$ taken into account in the previous section is now replaced by the interval $\left[e^{L}_{pq},e^{R}_{pq}\right]$ with $e^{L}_{pq}\leqslant e^{R}_{pq}$ and $e^{L}_{pq}=e^{R}_{pq}=e_{pq}$ if the information regarding the pair of levels $(p,q)\in P$ is precise. 

Of course, the starting point will be the definition of the \textit{interval consistency condition} taking into account the interval evaluations provided by the DM. 

\begin{condition}[Interval consistency condition]\label{Interval Consist Cond 1}
An interval pairwise comparison table is said consistent iff for all $(p,q)\in P$ there exists $e_{pq}\in\left[e^{L}_{pq},e^{R}_{pq}\right]$ such that the precise consistency Condition \ref{Precies Consistency Condition} holds.
\end{condition}

In more practical terms, the interval consistency condition above is satisfied iff there exists a precise pairwise comparison table $C=[e_{pq}]_{(p,q)\in P}$ such that $e_{pq}\in\left[e^{L}_{pq},e^{R}_{pq}\right]$ for all $(p,q)\in P$ and the consistency Condition \ref{Precies Consistency Condition} is satisfied. 

For each interval evaluation $\left[e^{L}_{pq},e^{R}_{pq}\right]$ we consider the new interval evaluation $\left[\overline{e}^{L}_{pq},\overline{e}^{R}_{pq}\right]$, where: 

\begin{itemize}[label={--}]
\item $e^{L}_{pq}+\delta^{L+}_{pq}-\delta^{L-}_{pq}=\overline{e}^{L}_{pq},$ for all $(p,q)\in P$, is the new left bound of the interval evaluation, 
\item $e^{R}_{pq}+\delta^{R+}_{pq}-\delta^{R-}_{pq}=\overline{e}^{R}_{pq},$ for all $(p,q)\in P$, is the new right bound of the interval evaluation,
\item $\delta^{L+}_{pq},\delta^{L-}_{pq},\delta^{R+}_{pq},\delta^{R-}_{pq}\in\Bbb{N}_{0}$ are auxiliary variables used to quantify how much the left or right bounds of the interval evaluation provided by the DM have to be modified. 
\end{itemize}

To check if the comparisons provided by the DM are consistent, we have to solve the following MILP problem

\begin{subequations}
\begin{gather}
\hspace{-5cm} y^{*}=\argmin z(y)=\displaystyle\sum_{(p,q) \; \in \; P}y_{pq}  \nonumber \\
\intertext{\hspace{3.75cm}subject to:}
\hspace{2,5cm}\overline{e}_{pk}+\overline{e}_{kq}+1=\overline{e}_{pq},  \qquad  (p,k), (k,q), (p,q) \in P,  \label{Interval Constraint_1} \\
\hspace{1cm}   e^{L}_{pq}+\delta^{L+}_{pq}-\delta^{L-}_{pq}=\overline{e}^{L}_{pq},  \qquad  (p,q)  \in P,  \label{Interval Constraint_1} \\
 \hspace{1cm}   e^{R}_{pq}+\delta^{R+}_{pq}-\delta^{R-}_{pq}=\overline{e}^{R}_{pq},  \qquad  (p,q)  \in P,  \label{Interval Constraint_2} \\
 \overline{e}^{L}_{pq}\leqslant\overline{e}_{pq}\leqslant\overline{e}^{R}_{pq}, \qquad (p,q)\in P, \label{Interval Constraint_1} \\
 \hspace{2,6cm}   \delta^{L+}_{pq}+\delta^{L-}_{pq}+\delta^{R+}_{pq}+\delta^{R-}_{pq}\leqslant  My_{pq},\;\;\;\; \; (p,q)  \in P,  \label{Interval Constraint_4} \\
  \hspace{3,5cm}  \delta_{pq}^{L-},\, \delta_{pq}^{L+},\, \delta_{pq}^{R-},\, \delta_{pq}^{R+} \in \Bbb{N}_0, \;\;\; \qquad (p,q)  \in P, \label{Interval Constraint_6} \\
   \hspace{1.5cm}  y_{pq}\in\{0,1\},\;\;\;\; \qquad (p,q)  \in P, \label{Interval Constraint_7} 
\end{gather}
\end{subequations}

\noindent where $M$ is a big positive number, while $y_{pq}$ is a binary variable used to check if some of the bounds of the interval $\left[e^{L}_{pq},e^{R}_{pq}\right]$ have to be modified or not. In particular, if $z^{*}=0$, then the information provided by the DM is consistent while, in the opposite case, some interval evaluation needs to be modified. Moreover, if $y_{pq}^{*}=0,$ then the bounds of the interval related to the pair of levels $(l_p,l_q)$ have not to be modified and the values $\overline{e}_{pq}\in\left[e^{L}_{pq},e^{R}_{pq}\right]$ will form a precise comparison table inferred by the imprecise one provided by the DM. If, instead, $y_{pq}^{*}=1$, the left or right bounds of the interval $\left[e^{L}_{pq},e^{R}_{pq}\right]$ need to be modified. In the following, we shall denote the MILP problem above with $MILP-I$.

If there exists a comparison table compatible with the information provided by the DM, it is reasonable asking if there is another possible comparison table or if the found comparison table is unique. For this reason, denoting by $\overline{e}^*_{pq}$ the optimal values obtained by solving $MILP-I$, we have to solve again the same problem with the addition of the following constraints 

\begin{subequations}
\begin{gather}
 \hspace{0.5cm}  \overline{e}_{pq}\leqslant \left(\overline{e}^{*}_{pq}-1\right)+My^{1}_{pq},  \qquad  (p,q)\in P, \label{eq: Constraint_1_Further_Int} \\ 
 \hspace{0.5cm}  \overline{e}_{pq}\geqslant \left(\overline{e}^{*}_{pq}+1\right)-My^{2}_{pq},  \qquad  (p,q)\in P, \label{eq: Constraint_2_Further_Int} \\ 
 \hspace{0.5cm}   \displaystyle\sum_{(p,q)\in P}\left(y^{1}_{pq}+y^{2}_{pq}\right)\leqslant 2\vert P\vert-1,\;\;\;\; \;   \label{eq: Constraint_3_Further_Int} \\
   \hspace{1.5cm}  y^{1}_{pq},y^{2}_{pq}\in\{0,1\},\;\;\;\; \qquad (p,q)  \in P\setminus P^{DM}, \label{eq: Constraint_4_Further_Int} 
\end{gather}
\end{subequations}

\noindent where constraints (\ref{eq: Constraint_1_Further_Int}) and  (\ref{eq: Constraint_2_Further_Int}) avoid that the MILP problem gives back the same $\overline{e}_{pq}^{*}$ value; if $y^{*}=0$, then there exists another precise comparison table that can be obtained from the interval evaluations provided by the DM. If, instead, $y^{*}>0$, then the precise comparison table found by $MILP-I$ is the unique one. Following the same procedure, in an iterative way, one can obtain all the precise comparison tables compatible with the imprecise information provided by the DM. Let us observe that, however, the number of precise comparison tables that can be extracted is equal to $\displaystyle\prod_{(p,q)\in P}\left(e_{pq}^{R}-e_{pq}^{L}+1\right)$ and, maybe, not all of them are consistent with respect to the consistency Condition \ref{Precies Consistency Condition}. If there exists more than one precise comparison table compatible with the information provided by the DM, then we can compute multicriteria preferences for each of them. In this case, taking into account robustness concerns, one could compute the frequency with which an alternative attains a certain rank position and the frequency with which an alternative is preferred to another. This probabilistic information can be seen as an application of the Stochastic Multicriteria Acceptability Analysis (SMAA) \citep{LahdelmaHokkanenSalminen1998,PelissariEtAl2019,TervonenFigueira2008} to our approach. 

Let us observe that if the DM provides imprecise information regarding all pairs of consecutive levels and we assume that the number of cards can be any positive number, not necessarily integer, then the preference information given by the DM is for sure precise and the procedure used to sample several sets of values compatible with this information is exactly the same introduced in \cite{CorrenteEtAl2017b}.

%%%%%%%%%%%%%%%%%%%%%%%%%%%%%%%%%%%%%%%
\subsection{An illustrative example}%%%
%%%%%%%%%%%%%%%%%%%%%%%%%%%%%%%%%%%%%%%
\noindent Let us suppose that the DM is able to provide the imprecise preference information shown in the comparison table below.

\[
\left\vert
\begin{array}{c|ccccc}\hline
           & l_{1} & l_{2} & l_{3} & l_{4} & l_{5} \\ \hline
    l_{1}  &\blacksquare & \left[0,1\right] & \left[5,6\right] & \left[5,8\right] & \left[7,11\right] \\
    l_{2}  & &\blacksquare & \left[1,2\right] & \left[4,6\right]  & \left[8,11\right]   \\
    l_{3}  & & &\blacksquare & \left[2,3\right]  & \left[6,8\right]      \\
    l_{4}  & & & &\blacksquare & \left[3,4\right]       \\
    l_{5}  & & & & &\blacksquare         \\  \hline
\end{array}
\right\vert
\]

Solving $MILP-{I}$, we obtain $z^{*}=1$ and $y^{*}_{13}=1$. In particular, we get $\delta^{L-*}_{13}=1$ meaning that, to restore the interval consistency, the interval $\left[5,6\right]$ should be replaced by the interval $\left[4,6\right]$ in which the old left bound has been decreased of one unit. Indeed, looking at the interval evaluations regarding pairs of levels $(l_1,l_2)$ and $(l_2,l_3)$, that are $\left[0,1\right]$ and $\left[1,2\right]$, one can easily observe that the maximum value for $\overline{e}_{13}$ should be equal to 4 (it is obtained when $\overline{e}_{12}=1$ and $\overline{e}_{23}=2$) but the minimum value for the comparison between levels $l_{1}$ and $l_{3}$ is, instead, $e^{L}_{13}=5$. Consequently, an interval pairwise comparison table consistent with the information provided by the DM is the following:

\[
\left\vert
\begin{array}{c|ccccc}\hline
           & l_{1} & l_{2} & l_{3} & l_{4} & l_{5} \\ \hline
    l_{1}  &\blacksquare & \left[0,1\right] & \left[\textcolor{red}{4},6\right] & \left[5,8\right] & \left[7,11\right] \\
    l_{2}  & &\blacksquare & \left[1,2\right] & \left[4,6\right]  & \left[8,11\right]   \\
    l_{3}  & & &\blacksquare & \left[2,3\right]  & \left[6,8\right]      \\
    l_{4}  & & & &\blacksquare & \left[3,4\right]       \\
    l_{5}  & & & & &\blacksquare         \\  \hline
\end{array}
\right\vert
\]

In particular, the solution of the previous $MILP-I$ problem will give back a precise consistent comparison table that can be extracted from the interval pairwise comparison table above, being composed of the values $\overline{e}^{*}_{pq}$ with $(p,q)\in P$. In this case, this precise comparison table is provided below: 

\[
\left\vert
\begin{array}{c|ccccc}\hline
           & l_{1} & l_{2} & l_{3} & l_{4} & l_{5} \\ \hline
    l_{1}  &\blacksquare & 1 & 4 & 7 & 11 \\
    l_{2}  & &\blacksquare & 2 & 5  & 9   \\
    l_{3}  & & &\blacksquare & 2  & 6      \\
    l_{4}  & & & &\blacksquare & 3       \\
    l_{5}  & & & & &\blacksquare         \\  \hline
\end{array}
\right\vert
\]

After replacing $e_{13}^{L}=5$ with $\overline{e}_{13}^{L}=4$, to check if it is possible extracting another precise consistent comparison table, we apply the same procedure described in the previous section solving the same problem with the addition of the constraints necessary to avoid to obtain the same precise comparison table shown above. In this case, $z^{*}=1$. This means that the precise comparison table obtained by solving $MILP-I$ is the unique one.

In this second part of the section, let us assume that the DM provides the following interval pairwise comparison table: 

\[
\left\vert
\begin{array}{c|ccccc}\hline
           & l_{1} & l_{2} & l_{3} & l_{4} & l_{5} \\ \hline
    l_{1}  &\blacksquare & \left[0,1\right] & \left[2,4\right] & \left[5,8\right] & \left[9,11\right] \\
    l_{2}  & &\blacksquare & \left[1,2\right] & \left[4,6\right]  & \left[8,11\right]   \\
    l_{3}  & & &\blacksquare & \left[2,3\right]  & \left[6,8\right]      \\
    l_{4}  & & & &\blacksquare & \left[3,4\right]       \\
    l_{5}  & & & & &\blacksquare         \\  \hline
\end{array}
\right\vert
\]

Solving the $MILP-I$ problem, we obtain $z^{*}=0$ meaning that it is possible extracting at least one consistent precise comparison table and, in particular, the precise comparison table composed of the values $\overline{e}_{pq}$ obtained by the $MILP-I$ problem is the following: 

\[
\left\vert
\begin{array}{c|ccccc}\hline
           & l_{1} & l_{2} & l_{3} & l_{4} & l_{5} \\ \hline
    l_{1}  &\blacksquare & 0 & 2 & 5 & 9 \\
    l_{2}  & &\blacksquare & 1 & 4 & 8   \\
    l_{3}  & & &\blacksquare & 2  & 6      \\
    l_{4}  & & & &\blacksquare & 3       \\
    l_{5}  & & & & &\blacksquare         \\  \hline
\end{array}
\right\vert
\]

By avoiding to obtain again the same precise consistent comparison table, we are able to find other 10 precise comparison tables\footnote{Let us observe that the number of precise comparison tables that can be extracted from the imprecise comparison table provided by the DM are 20,736.} that can be extracted from the imprecise comparison table provided by the DM and that are consistent. For saving space we will not include all of them in the paper but it is possible to download the file containing them clicking on the following link: \href{http://www.antoniocorrente.it/wwwsn/images/allegati_articoli/SupplementaryDCM.zip}{further consistent comparison tables}.

%%%%%%%%%%%%%%%%%%%%%%%%%%%%%%%%%%%%%%%%%%%%%%%%%%%%%%%%%%%%%%%%%%%%%%%%%%
\subsection{Imprecise and missing information}\label{Missing Imprecise}%%%
%%%%%%%%%%%%%%%%%%%%%%%%%%%%%%%%%%%%%%%%%%%%%%%%%%%%%%%%%%%%%%%%%%%%%%%%%%
\noindent If the DM is not sure about the comparison between different levels at hand, as observed in the previous section, he can decide to provide some imprecise evaluations expressed as intervals or, in the extreme case, he can avoid to provide any preference information regarding some pairs of levels. For this reason, in this section we shall consider the case in which the information provided by the DM is expressed as interval of possible values or as missing values. 

Using the notation introduced in the previous sections, we shall denote by $P^{DM}$ the set composed of the pairs of levels for which the DM provided an interval comparison, while for the remaining pairs of levels (those in $P\setminus P^{DM}$), the information is missing. The objective will be therefore, finding values $\overline{e}_{pq}\in\left[e_{pq}^{L},e_{pq}^{R}\right]$ for all $(p,q)\in P^{DM}$ and $\overline{e}_{pq}\in\Bbb{N}_0$ for $(p,q)\in P\setminus P^{DM}$ such that the precise consistency Condition \ref{Precies Consistency Condition} is fulfilled. 

To check if there exists a precise comparison table that can be extracted from the imprecise and missing information provided by the DM, one has to solve the following MILP problem that will be denoted by $MILP-MI$:

\begin{subequations}
\begin{gather}
\hspace{-5cm} y^{*}=\argmin z(y)=\displaystyle\sum_{(p,q) \; \in \; P^{DM}}y_{pq}  \nonumber \\
\intertext{\hspace{3.75cm}subject to:}
\hspace{1.9cm}\overline{e}_{pk}+\overline{e}_{kq}+1=\overline{e}_{pq},  \qquad  (p,k), (k,q), (p,q) \in P,  \label{Interval Miss Constraint_1} \\
\hspace{1cm}   e^{L}_{pq}+\delta^{L+}_{pq}-\delta^{L-}_{pq}=\overline{e}^{L}_{pq},  \qquad  (p,q)  \in P^{DM},  \label{Interval Miss Constraint_2} \\
 \hspace{1cm}   e^{R}_{pq}+\delta^{R+}_{pq}-\delta^{R-}_{pq}=\overline{e}^{R}_{pq},  \qquad  (p,q)  \in P^{DM},  \label{Interval Miss Constraint_3} \\
 \overline{e}^{L}_{pq}\leqslant\overline{e}_{pq}\leqslant\overline{e}^{R}_{pq}, \qquad (p,q)\in P^{DM}, \label{Interval Miss Constraint_4} \\
 \hspace{2,6cm}   \delta^{L+}_{pq}+\delta^{L-}_{pq}+\delta^{R+}_{pq}+\delta^{R-}_{pq}\leqslant  My_{pq},\;\;\;\; \; (p,q)  \in P^{DM},  \label{Interval Miss Constraint_5} \\
  \hspace{3,5cm}  \delta_{pq}^{L-},\, \delta_{pq}^{L+},\, \delta_{pq}^{R-},\, \delta_{pq}^{R+} \in \Bbb{N}_0, \;\;\; \qquad (p,q)  \in P^{DM}, \label{Interval Miss  Constraint_6} \\
  \hspace{1.5cm}  \overline{e}_{pq}\in \Bbb{N}_0, \;\;\; \qquad (p,q)  \in P, \label{Interval Miss  Constraint_7} \\
   \hspace{2.5cm}  y_{pq}\in\{0,1\}, \;\;\qquad (p,q)  \in P^{DM}. \label{Interval Constraint_7} 
\end{gather}
\end{subequations}

If $z^{*}=0$, then it is possible to find values $\overline{e}_{pq}\in\left[e_{pq}^{L},e_{pq}^{R}\right]$ for all $(p,q)\in P^{DM}$ and values $\overline{e}_{pq}\in \Bbb{N}_0$ for all $(p,q)\in P\setminus P^{DM}$ so that the comparison table is consistent with respect to Condition \ref{Precies Consistency Condition}. In the opposite case ($z^{*}>0$) some interval information provided by the DM has to be modified to make the comparison table consistent. If there exists at least one precise comparison table that can be extracted from the imprecise and missing information, to check for others, one has to solve the same problem with the addition of the constraints (\ref{eq: Constraint_1_Further_Int})-(\ref{eq: Constraint_4_Further_Int}) introduced in section \ref{sec:Imprecise}. A null $z^{*}$ value implies that there exists another precise comparison table compatible with the imprecise and missing information provided by the DM. In the opposite case ($z^{*}>0$), the previously found precise comparison table is the unique one.

%%%%%%%%%%%%%%%%%%%%%%%%%%%%%%%%%%%%%%%
\subsection{An illustrative example}%%%
%%%%%%%%%%%%%%%%%%%%%%%%%%%%%%%%%%%%%%%
\noindent Let us assume that the DM was able to provide the following imprecise and missing comparison table: 

\[
\left\vert
\begin{array}{c|ccccc}\hline
           & l_{1} & l_{2} & l_{3} & l_{4} & l_{5} \\ \hline
    l_{1}  &\blacksquare & ? & \left[2,4\right] & \left[5,8\right] & \left[7,11\right] \\
    l_{2}  & &\blacksquare & \left[1,2\right] & ?  & \left[6,9\right]   \\
    l_{3}  & & &\blacksquare & \left[2,3\right]  & ?      \\
    l_{4}  & & & &\blacksquare & \left[1,4\right]       \\
    l_{5}  & & & & &\blacksquare         \\  \hline
\end{array}
\right\vert
\]

By solving the previous $MILP-MI$ problem, we find that it is feasible and the optimal objective value is $z^{*}=0$. This means that the information provided by the DM is consistent and that there exists at least one consistent precise comparison table. The precise comparison table obtained solving the $MILP-MI$ problem is shown below.

\[
\left\vert
\begin{array}{c|ccccc}\hline
           & l_{1} & l_{2} & l_{3} & l_{4} & l_{5} \\ \hline
    l_{1}  &\blacksquare & 0 & 2 & 5 & 7 \\
    l_{2}  & &\blacksquare & 1 & 4  & 6   \\
    l_{3}  & & &\blacksquare & 2  & 4      \\
    l_{4}  & & & &\blacksquare & 1       \\
    l_{5}  & & & & &\blacksquare         \\  \hline
\end{array}
\right\vert
\]

To avoid that the same solution is obtained, we have added constraints (\ref{eq: Constraint_1_Further_Int})-(\ref{eq: Constraint_4_Further_Int}) obtaining 28 different precise comparison tables compatible with the imprecise and missing preference information provided by the DM. Again, we do not include all of them in the manuscript but the interested reader can download them clicking on the following link: \href{http://www.antoniocorrente.it/wwwsn/images/allegati_articoli/SupplementaryDCM.zip}{further consistent comparison tables}.

%%%%%%%%%%%%%%%%%%%%%
\section{The approach of the comparison table based on the DCM for the Choquet integral preference model: an illustrative example}\label{sec:MAVT}
%%%%%%%%%%%%%%%%%%%%%
\noindent We consider the example in \cite{BotteroEtAl18} regarding the selection of the best strategy for the requalification of an abandoned quarry. The set of criteria is as follows:

\begin{enumerate}
  \item \textit{Investment costs} (Scale unit: $K$\euro; Code: {\tt COSTS}; notation: $g_1$; preference direction: minimization). This criterion comprises the requalification costs of the quarry.
  \item \textit{Profitability} (Scale unit: verbal levels (seven); Code: {\tt PROFI}; notation: $g_2$; preference direction: minimization). This criterion comprises the future income the project is expected to produce for the local population.
  \item \textit{New services for the population} (Scale unit: verbal levels (seven); Code: {\tt SERVI}; notation: $g_3$; preference direction: maximization). This criterion models the possibility of recruiting workers.
  \item \textit{Naturalized surface} (Scale unit: hectares; Code: {\tt SURFA}; notation: $g_4$; preference direction: maximization). This criterion comprises the impacts of a project in the landscape quality and bio-diversity conservation.
  \item \textit{Environmental effects} (Scale unit: verbal levels (seven); Code: {\tt ENVIR}; notation: $g_5$; preference direction: maximization). This criterion comprises the impacts of a project in the environmental system.
   \item \textit{Consistency with local planning requirements} (Scale unit: two levels (Yes-No); Code: {\tt CONSI}; notation: $g_6$; preference direction: maximization). This criterion is related to the administrative feasibility of the project with respect to some urban constraints (if it is feasible, the answer is ``yes'' [1]; otherwise, the answer is ``no'' [0].
\end{enumerate}

The verbal scale used for the criteria $g_2$, $g_3$, and $g_5$ comprises the following seven levels  (in between parenthesis we used a numerical code for each level): {\tt very bad}[1]; {\tt bad}[2]; {\tt rather bad}[3]; {\tt average}[4]; {\tt rather good}[5]; {\tt good}[6]; {\tt very good}[7].

There are five potential requalification projects. They can be summarized as follows: basic reclamation ($a_1$); valuable forest($a_2$); wetland ($a_3$); ecological network ($a_4$); and, multi-functional area ($a_5$). The performance table can be presented as follows.

\begin{table}[htb!]
  \centering
    {\small
    \begin{tabular}{ccccccc}\hline
$a$ & {\tt COSTS} ($g_1(a)$) & {\tt PROFI} ($g_2(a)$) & {\tt SERVI} ($g_3(a)$) & {\tt SURFA} ($g_4(a)$) & {\tt ENVIR} ($g_5(a)$)  & {\tt CONSI} ($g_6(a)$) \\ \hline
  $a_1$ & $\;\,30\;000$ & 3 & 1 & 2.0 & 4 & 1 \\
  $a_2$ & $\;\,45\;000$ & 3 & 5 & 5.0 & 5 & 1 \\
  $a_3$ & $\;\,90\;000$ & 1 & 6 & 3.2 & 7 & 1 \\
  $a_4$ & $120\;000$    & 1 & 7 & 3.5 & 6 & 1 \\
  $a_5$ & $900\;000$    & 7 & 7 & 1.0 & 3 & 0 \\ \hline
    \end{tabular}
    }
  \caption{Performance table}\label{tab:Perf_Table}
\end{table}

%%%%%%%%%%%%%%%%%%%%%%%%%%%%%%%%%%%%%%%%%%%%%%%%%%%%%%%%%%%%%%
\subsection{The multicriteria evaluation of the projects}\label{Utilities}%%%
%%%%%%%%%%%%%%%%%%%%%%%%%%%%%%%%%%%%%%%%%%%%%%%%%%%%%%%%%%%%%%
\noindent In this section, we shall apply the extension of the DCM approach described in the previous sections to each criterion, to build a unique common scale $[0,100]$ and, then, aggregating these utility values by means of the Choquet integral preference model.

\begin{enumerate}
    \item Criterion $g_1$  {\tt COSTS} (Investment costs).
    For this criterion the DM considered as the lowest cost $0$ and the highest $1000$ providing the following table regarding comparisons between different levels. 
    
    \[
    \left\vert
    \begin{array}{c|ccccc}\hline
        & l_{1}(1000)  & l_{2}(750) & l_{3}(500) & l_{4}(250) & l_{5}(0) \\ \hline
        l_{1}(1000) & \blacksquare & 3 & 6 & 8 & 9 \\
        l_{2}(550)~ & & \blacksquare   & 2 & 4 & 5 \\
        l_{3}(500)~ & &   &  \blacksquare  & 1 & 2 \\
        l_{4}(250)~ & &   &   & \blacksquare   & 0 \\
        l_{5}(0) & &   &   &   &  \blacksquare \\ \hline
        \end{array}
        \right\vert
    \]
    
    As one can easily check, the table satisfies the precise consistency Condition \ref{Precies Consistency Condition}. This means that there is not the necessity to revise the information provided by the DM and, consequently, we can proceed to assign a single value to the considered levels. Let us proceed by computing the number of units between the lowest level and the highest one, $h=e_{15}+1=10$ and, therefore, $u(l_{5})=u(l_{1})+10\cdot \alpha$. Fixing the reference levels $u_1(l_{1}=1000)=0$ and $u_1(l_{5}=0)=100$ the value of the unit is $\alpha = \frac{u(l_{5})-u(l_{1})}{10}=\frac{100-0}{10} = 10$. The remaining values are computed as follows:
    \begin{itemize}
      \item[~]  $u_1(l_2=750)= u(l_{1})+(e_{12}+1)\cdot \alpha = 0 + (3+1)\cdot 10 = 40$.
      \item[~]  $u_1(l_3=500)= u(l_{1})+(e_{13}+1)\cdot \alpha = 0 + (6+1)\cdot 10 = 70$.
      \item[~]  $u_1(l_4=250)= u(l_{1})+(e_{14}+1)\cdot \alpha = 0 + (8+1)\cdot 10 = 90$.
    \end{itemize}
    In order to compute the values of the performances we need to proceed with a linear interpolation.
     \begin{itemize}
      \item[~]  $u_1(30)= 100 + ((30 - 0)/(250-0))(90-100) = 98.8$.
      \item[~]  $u_1(45)= 100 + ((45-0)/(250-0))(90-100)= 98.2$.
      \item[~]  $u_1(90)= 100 + ((90-0)/(250-0))(90-100)=96.4$.
      \item[~]  $u_1(120)= 100 + ((120-0)/(250-0))(90-100)= 95.2$.
      \item[~]  $u_1(900)= 40 + (900-750)/(1000-750))(0-40) = 16$.
    \end{itemize}
    These values can now be inserted in the scores table.
    \item Criterion $g_2$  {\tt PROFI} (Profitability). With respect to criterion $g_{2}$, the DM provided the imprecise and missing pairwise comparison table shown below. 

    \[
    \left\vert
    \begin{array}{c|ccccccc}\hline
                  & l_{1}(vb) & l_{2}(b) & l_{3}(rb) & l_{4}(a) & l_{5}(rg) & l_{6}(g) & l_{7}(vg) \\ \hline
        l_{1}(vb) &  \blacksquare & \left[0,1\right] & 2 & 4 & 8 & 10 & 14 \\
        l_{2}(b)~ &  & \blacksquare   & 1 & \left[3,4\right] & 6 & 9 & 13 \\
        l_{3}(rb) &  &   &  \blacksquare  & 1 & 4 & \left[6,7\right] & 11 \\
        l_{4}(a)~ &  &   &   & \blacksquare   & ? & 5 & 9 \\
        l_{5}(rg) &  &   &   &   & \blacksquare   & 2 & ? \\
        l_{6}(g)~ &  &   &   &   &   & \blacksquare   & 3 \\
        l_{7}(vg) &  &   &   &   &   &    & \blacksquare  \\ \hline
        \end{array}
        \right\vert
    \]
    Solving the $MILP-MI$, we find $z^{*}=1$. In particular, $y_{15}^{*}=1$ and $\delta_{15}^{L-*}=1$. This means that the comparison table provided by the DM is inconsistent and that the consistency can be restored by modifying the comparison between levels $l_{1}$ and $l_{5}$ from $8$ to the interval $\left[7,8\right]$. Indeed, for the consistency Condition \ref{Precies Consistency Condition}, it should hold that $e_{13}+e_{35}+1=e_{15}$. Anyway, in the comparison table provided by the DM this is not true since $e_{13}+e_{35}+1=7$, while $e_{15}=8$. The DM agreed to modify this comparison and, consequently, we checked again if the comparison table is consistent. In this case, $z^{*}=0$ and, therefore, there exists at least one precise comparison table compatible with the provided information that is shown below: 

    \[
    \left\vert
    \begin{array}{c|ccccccc}\hline
                  & l_{1}(vb) & l_{2}(b) & l_{3}(rb) & l_{4}(a) & l_{5}(rg) & l_{6}(g) & l_{7}(vg) \\ \hline
        l_{1}(vb) &  \blacksquare & 0 & 2 & 4 & 7 & 10 & 14 \\
        l_{2}(b)~ &  & \blacksquare   & 1 & 3 & 6 & 9 & 13 \\
        l_{3}(rb) &  &   &  \blacksquare  & 1 & 4 & 7 & 11 \\
        l_{4}(a)~ &  &   &   & \blacksquare   & 2 & 5 & 9 \\
        l_{5}(rg) &  &   &   &   & \blacksquare   & 2 & 6 \\
        l_{6}(g)~ &  &   &   &   &   & \blacksquare   & 3 \\
        l_{7}(vg) &  &   &   &   &   &    & \blacksquare  \\ \hline
        \end{array}
        \right\vert
    \]

    Checking for another precise comparison table, we found $z^{*}>0$. Therefore, the above table is the unique one.  As for the previous criterion, we fix $u_2(vb=1) = 0$ and $u_2(vg=7)=100$. Then, looking at the comparison table we have $h = e_{17}+1 = 15$ and $u_{2}(vg=7)=u_{2}(vb=1)+15\cdot \alpha$ and, consequently $\alpha = \frac{u_{2}(vg=7)-u_{2}(vb=1)}{15}=\frac{100-0}{15} = 6.667$.  The other five levels are computed as follows:
    \begin{itemize}
      \item[~]  $u_2(b=2)= u_2(vb=1) + (e_{12}+1)\cdot \alpha = 0 + 1\cdot 6.667 = 6.667$,
      \item[~]  $u_2(rb=3)= u_2(vb=1) + (e_{13}+1)\cdot \alpha = 0 + 3\cdot 6.667 = 20$,
      \item[~]  $u_2(a=4)= u_2(vb=1) + (e_{14}+1)\cdot \alpha = 0 + 5\cdot 6.667 = 33.333$,
      \item[~]  $u_2(rg=5)= u_2(vb=1) + (e_{15}+1)\cdot \alpha = 0 + 8\cdot 6.667 = 53.333$,
      \item[~]  $u_2(g=6)= u_2(vb=1) + (e_{16}+1)\cdot \alpha = 0 + 11\cdot 6.667 = 73.333$.
    \end{itemize}
    \item Criterion $g_3$  {\tt SERVI} (New services for the population). Also in this case, the DM provided imprecise preference information as shown in the following comparison table: 
     \[
    \left\vert
    \begin{array}{c|ccccccc}\hline
                  & l_{1}(vb) & l_{2}(b) & l_{3}(rb) & l_{4}(a) & l_{5}(rg) & l_{6}(g) & l_{7}(vg) \\ \hline
        l_{1}(vb) &  \blacksquare & \left[1,2\right] & \left[3,4\right] & \left[6,7\right] & \left[9,10\right] & \left[13,14\right] & \left[18,19\right] \\
        l_{2}(b)~ &  & \blacksquare   & \left[1,2\right] & \left[4,5\right] & \left[7,8\right] & \left[11,12\right] & \left[16,17\right] \\
        l_{3}(rb) &  &   &  \blacksquare  & \left[2,3\right] & \left[5,6\right] & \left[9,10\right] & \left[14,15\right] \\
        l_{4}(a)~ &  &   &   & \blacksquare   & \left[2,3\right] & \left[6,7\right] & \left[11,12\right] \\
        l_{5}(rg) &  &   &   &   & \blacksquare   & \left[3,4\right] & \left[8,9\right] \\
        l_{6}(g)~ &  &   &   &   &   & \blacksquare   & \left[4,5\right] \\
        l_{7}(vg) &  &   &   &   &   &    & \blacksquare  \\ \hline
        \end{array}
        \right\vert
    \]
    Solving the $MILP-I$ problem, we find $z^{*}=0$. This means that the information provided by the DM is intervally consistent and there exists at least one precise comparison table that can be extracted from the imprecise one. The precise comparison table obtained by solving the $MILP-I$ problem is shown below: 

     \[
    \left\vert
    \begin{array}{c|ccccccc}\hline
                  & l_{1}(vb) & l_{2}(b) & l_{3}(rb) & l_{4}(a) & l_{5}(rg) & l_{6}(g) & l_{7}(vg) \\ \hline
        l_{1}(vb) &  \blacksquare & 1 & 3 & 6 & 9 & 13 & 18 \\
        l_{2}(b)~ &  & \blacksquare   & 1 & 4 & 7 & 11 & 16 \\
        l_{3}(rb) &  &   &  \blacksquare  & 2 & 5 & 9 & 14 \\
        l_{4}(a)~ &  &   &   & \blacksquare   & 2 & 6 & 11 \\
        l_{5}(rg) &  &   &   &   & \blacksquare   & 3 & 8 \\
        l_{6}(g)~ &  &   &   &   &   & \blacksquare   & 4 \\
        l_{7}(vg) &  &   &   &   &   &    & \blacksquare  \\ \hline
        \end{array}
        \right\vert
    \]

    This time, checking for other possible precise comparison tables compatible with the imprecise information provided by the DM we find 6 others. For space reasons, we do not report all of them in the manuscript but we shall take into account all these tables to obtain the final recommendation on the considered problem. These precise tables can be found clicking on the following link: \href{http://www.antoniocorrente.it/wwwsn/images/allegati_articoli/SupplementaryDCM.zip}{further consistent comparison tables}. To obtain the utilities of the 7 levels on $g_{3}$, we consider the precise table above. We fix $u_3(vb=1) = 0$ and $u_3(vg=7)=100$. Then, compute $h = e_{17}+1 = 19$ and $u_{3}(vg=7)=u_{3}(vb=1)+19\cdot \alpha$. Consequently, the value of the unit $\alpha = \frac{u_{3}(vg=7)-u_{3}(vb=1)}{19}=(100-0)/19 = 5.263$.  The other five levels are computed as follows:
    \begin{itemize}
      \item[~]  $u_3(b=2)= u_{3}(vb=1) + (e_{12}+1)\cdot \alpha= 0 + 2\cdot 5.263 = 10.526$,
      \item[~]  $u_3(rb=3)= u_{3}(vb=1) + (e_{13}+1)\cdot \alpha = 0 + 4\cdot 5.263 = 21.052$,
      \item[~]  $u_3(a=4)= u_{3}(vb=1) + (e_{14}+1)\cdot \alpha = 0 + 7\cdot 5.263 = 36.841$,
      \item[~]  $u_3(rg=5)= u_{3}(vb=1) + (e_{15}+1)\cdot \alpha = 0 + 10\cdot 5.263 = 52.630$,
      \item[~]  $u_3(g=6)= u_{3}(vb=1) + (e_{16}+1)\cdot \alpha = 0 + 14\cdot 5.263 = 73.682$.
    \end{itemize}
    \item Criterion $g_4$  {\tt SURFA} (Naturalized surface). For this criterion we consider the following comparison table and $u_4(l_1=1.0)=0$ and $u_4(l_5=5.0)=100$.
    \[
    \left\vert
    \begin{array}{c|ccccc}\hline
        & l_{1}(1.0) & l_{2}(2.9) & l_{3}(3.2) & l_{4}(3.5) & l_{5}(5.0) \\ \hline
        l_{1}(1.0) & \blacksquare & 1 & 3 & 5 & 8 \\
        l_{2}(2.0) & & \blacksquare   & 1 & 3 & 6 \\
        l_{3}(3.2) & &   &  \blacksquare  & 1 & 4 \\
        l_{4}(3.5) & &   &   & \blacksquare   & 2 \\
        l_{5}(5.0) & &   &   &   &  \blacksquare \\ \hline
        \end{array}
        \right\vert
    \]
    The comparison table is consistent and, therefore, we can compute the utilities of the different levels. The number of units between the two extreme levels is $h=e_{15}+1=9$, while $u_{4}(l_{5}=5.0)=u_{4}(l_{1}=1.0)+9\cdot \alpha$. Consequently, the value of the unit is $\alpha=\frac{u_{4}(l_{5}=5.0)-u_{4}(l_{1}=1.0)}{9}=\frac{100-0}{9} = 11.111$. With this information we can easily calculate the value of the remaining levels:
    \begin{itemize}
      \item[~]  $u_4(l_2=2.9)= u_{4}(l_{1}=1.0) + (e_{12}+1)\cdot \alpha = 0 + 2\cdot 11.111 = 22.222$,
      \item[~]  $u_4(l_3=3.2)= u_{4}(l_{1}=1.0) + (e_{13}+1)\cdot \alpha = 0 + 4\cdot 11.111 = 44.444$,
      \item[~]  $u_4(l_4=3.5)= u_{4}(l_{1}=1.0) + (e_{14}+1)\cdot \alpha = 0 + 6\cdot 11.111 = 66.666$.
    \end{itemize}
    These values can now be inserted in the scores table.
    \item Criterion $g_5$  {\tt ENVIR} (Environmental effects). Regarding this criterion, the DM provides interval and missing information as shown in the comparison table below:
    \[
    \left\vert
    \begin{array}{c|ccccccc}\hline
                  & l_{1}(vb) & l_{2}(b) & l_{3}(rb) & l_{4}(a) & l_{5}(rg) & l_{6}(g) & l_{7}(vg) \\ \hline
        l_{1}(vb) &  \blacksquare & ? & \left[0,1\right] & \left[2,3\right] & \left[4,5\right] & \left[6,7\right] & \left[9,10\right] \\
        l_{2}(b)~ &  & \blacksquare   & ? & \left[1,2\right] & \left[3,4\right] & \left[5,6\right] & \left[8,9\right] \\
        l_{3}(rb) &  &   &  \blacksquare & ? & \left[2,3\right] & \left[4,5\right] & \left[7,8\right]  \\
        l_{4}(a)~ &  &   &   & \blacksquare   & ? & \left[2,3\right] & \left[5,6\right] \\
        l_{5}(rg) &  &   &   &   & \blacksquare   & ? & \left[3,4\right] \\
        l_{6}(g)~ &  &   &   &   &   & \blacksquare   & ? \\
        l_{7}(vg) &  &   &   &   &   &    & \blacksquare  \\ \hline
        \end{array}
        \right\vert
    \]

    Solving the $MILP-MI$ problem, we find $z^{*}=0$. Consequently, there exists at least one precise comparison table compatible with the information provided by the DM. The precise table obtained solving the $MILP-MI$ problem is the following: 

    \[
    \left\vert
    \begin{array}{c|ccccccc}\hline
                  & l_{1}(vb) & l_{2}(b) & l_{3}(rb) & l_{4}(a) & l_{5}(rg) & l_{6}(g) & l_{7}(vg) \\ \hline
        l_{1}(vb) &  \blacksquare & 0 & 1 & 3 & 5 & 7 & 10 \\
        l_{2}(b)~ &  & \blacksquare   & 0 & 2 & 4 & 6 & 9 \\
        l_{3}(rb) &  &   &  \blacksquare & 1 & 3 & 5 & 8  \\
        l_{4}(a)~ &  &   &   & \blacksquare   & 1 & 3 & 6 \\
        l_{5}(rg) &  &   &   &   & \blacksquare   & 1 & 4 \\
        l_{6}(g)~ &  &   &   &   &   & \blacksquare   & 2 \\
        l_{7}(vg) &  &   &   &   &   &    & \blacksquare  \\ \hline
        \end{array}
        \right\vert
    \]
    
    Also in this case, we found other 7 precise comparison tables compatible with the provided information. We do not report them here but we shall take into account all of them to provide the final recommendation. All tables can be consulted clicking on the following link: \href{http://www.antoniocorrente.it/wwwsn/images/allegati_articoli/SupplementaryDCM.zip}{further consistent comparison tables}. We show how utilities for this criterion are computed taking into account the precise comparison table obtained by solving the $MILP-MI$ problem.  We proceed in the same way as for criteria $g_2$ and $g_3$. Therefore, $u_5(vb=1) = 0$, $u_5(vg=7)=100$ and $h = e_{17}+1 = 11$. Consequently, $u_{5}(vg=7) = u_{5}(vb=1) + 11\cdot \alpha$ and $\alpha = \frac{u_{5}(vg=7)-u_{5}(vb=1)}{11} = (100-0)/11 = 9.091$.  The other five levels are computed as follows:
    \begin{itemize}
      \item[~]  $u_5(b=2)= u_{5}(vb=1) + (e_{12}+1)\cdot \alpha = 0 + 1\cdot 9.091 = 9.091$,
      \item[~]  $u_5(rb=3)= u_{5}(vb=1) + (e_{13}+1)\cdot \alpha = 0 + 2\cdot 9.091 = 18.182$,
      \item[~]  $u_5(a=4)= u_{5}(vb=1) + (e_{14}+1)\cdot \alpha = 0 + 4\cdot 9.091 = 36.364$,
      \item[~]  $u_5(rg=5)= u_{5}(vb=1) + (e_{15}+1)\cdot \alpha = 0 + 6\cdot 9.091 = 54.546$,
      \item[~]  $u_5(g=6)= u_{5}(vb=1) + (e_{16}+1)\cdot \alpha = 0 + 8\cdot 9.091 = 72.728$.
    \end{itemize}
    \item Criterion $g_6$  {\tt CONSI} (Consistency with local requirements). For this criterion we simple consider that ``no'' has a zero value and ``yes'', a 100 value.
\end{enumerate}

The table containing the utilities obtained by the previous steps is shown below. Let us observe that for criteria $g_{3}$ and $g_{5}$ we reported the utilities obtained solving the corresponding $MILP-MI$ problems. 

\begin{table}[htb!]
  \centering
    \begin{tabular}{ccccccc}\hline
$a$ & $u_1(g_1(a))$ & $u_2(g_2(a))$ & $u_3(g_3(a))$ & $u_4(g_4(a))$ & $u_5(g_5(a))$  & $u_6(g_6(a))$ \\ \hline
  $a_1$ & 98.8 & 20.00 & 0.0000 & 11.111 & 36.364 & 100.00 \\
  $a_2$ & 98.2 & 20.00 & 52.630 & 100.00 & 54.546 & 100.00 \\
  $a_3$ & 96.4 & 0.000 & 73.682 & 22.222 & 100.00 & 100.00 \\
  $a_4$ & 95.2 & 0.000 & 100.00 & 66.666 & 72.728 & 100.00 \\
  $a_5$ & 16.0 & 100.0 & 100.00 & 0.000 & 18.182 & 0.0000 \\ \hline
    \end{tabular}
  \caption{Scores table}\label{tab:Perf_Table}
\end{table}

\subsection{The capacities of the Choquet integral}
\noindent As already explained in section \ref{sec:Ratio}, we shall use the DCM approach to get the capacity of the Choquet integral. For such a reason, we shall consider the following dummy projects: 
\begin{itemize}[label={--}]
    \item $p_{1}=(100,0,0,0,0,0)$, $p_{2}=(0,100,0,0,0,0)$, $p_{3}=(0,0,100,0,0,0)$, $p_{4}=(0,0,0,100,0,0)$, $p_{5}=(0,0,0,0,100,0)$, $p_{6}=(0,0,0,0,0,100)$,
    \item $p_{15}=(100,0,0,0,100,0)$ and $p_{45}=(0,0,0,100,100,0)$.
\end{itemize}

The DM is able to give the following ranking between the projects: 

$$
p_{6}\prec p_{1} \prec p_{5} \prec \left\{p_2,p_3\right\} \prec p_4 \prec p_{45} \prec p_{15}.
$$

This means that $p_{6}$ is the least important project, $p_{15}$ is the most important one and projects $p_{2}$ and $p_{3}$ have exactly the same importance. In terms of criteria, there is interaction between criteria $g_{4}$ and $g_{5}$ as well as between criteria $g_{1}$ and $g_{5}$. The pairwise comparison table provided by the DM is the following:

    \[
    \left\vert
    \begin{array}{c|ccccccc}\hline
                    & \{p_6\} & \{p_1\} & \{p_5\} & \{p_2,p_3\} & \{p_4\} & \{p_{45}\} & \{p_{15}\} \\ \hline
        \{p_6\}     &  \blacksquare & 1 &  3 & 4 & 6 & 9 & 14 \\
        \{p_1\}     &  & \blacksquare   & 1 &  2 & 4 & 7 & 12 \\
        \{p_5\}     &  &   &  \blacksquare  & 0  & 2 & 5 & 10 \\
        \{p_2,p_3\} &  &   &   & \blacksquare    & 1 & 4 &  9\\
        \{p_4\}     &  &   &   &   & \blacksquare   & 2 & 7 \\
        \{p_{45}\}  &  &   &   &   &   & \blacksquare   & 4 \\
        \{p_{15}\}  &  &   &   &   &   &    & \blacksquare   \\ \hline
        \end{array}
        \right\vert
    \]

    The comparisons provided by the DM are precise and the consistency Condition \ref{Precies Consistency Condition} is satisfied. Therefore, it is not necessary revising the given information. To get the ratio between the importance of the most important project $(p_{15})$ and the least important one $(p_{6})$, the DM is asked the following question: ``If you give a single vote to project $p_6$ how many can you assign to project $p_{15}$"? The DM says 8 and, consequently, project $p_{15}$ is 8 times more important than project $p_6$. \\
    Knowing the $z$ value, is therefore possible obtaining the values of the capacities as well as the M\"{o}bius parameters necessary to apply the Choquet integral preference model. One gets the values given in the table below:

   \[
    \left\vert
    \begin{array}{c|cc}\hline
                    & \mu(\cdot) & m(\cdot) \\ \hline
        \{p_6\}     & 0.0541 & 0.0541 \\
        \{p_1\}     & 0.1046 & 0.1046 \\
        \{p_5\}     & 0.1552 & 0.1552 \\
        \{p_2,p_3\} & 0.1805 & 0.1805 \\
        \{p_4\}     & 0.2310 & 0.2310 \\
        \{p_{45}\}  & 0.3068 & -0.0794  \\
        \{p_{15}\}  & 0.4332 & 0.1732 \\ \hline
        \end{array}
        \right\vert
    \]

Let us underline that from the values in the table above, on one hand, criteria $g_{4}$ and $g_{5}$ are negatively interacting, while, on the other hand, $g_{1}$ and $g_{5}$ are positively interacting. 

At this point, considering the scores obtained in the previous subsection and the capacities shown in the previous table, we find the ranking of the 5 projects as follows: 

$$
C_{\mu}(a_1) = 32.9999 < C_{\mu}(a_5) = 43.3712 < C_{\mu}(a_4) = 67.4332 < 
$$
$$
< C_{\mu}(a_3) = 63.8491 < C_{\mu}(a_4) = 67.4332 
$$

\noindent for which $a_{1}$ is the worst project, while $a_{4}$ is the best one. 

%%%%%%%%%%%%%%%%%%%%%%%%%%%%%%%%%%%%%%%%%%%%%%%%%%%%%%%%%%%%%%%%%
\subsection{On the use of SMAA for obtaining more robust conclusions}%%%%
%%%%%%%%%%%%%%%%%%%%%%%%%%%%%%%%%%%%%%%%%%%%%%%%%%%%%%%%%%%%%%%%%
\noindent As already observed in the previous section, regarding criteria $g_{3}$ and $g_{5}$, there exists more than one comparison table compatible with the imprecise or missing information provided by the DM. Comprehensively, we have 7 comparison tables for $g_{3}$, 8 comparison tables for $g_{5}$ and only one comparison table for the remaining criteria. Since using each comparison table we can obtain one utility for the considered levels, we need to take into account 56 different cases. All the matrices can be consulted clicking on the following link: \href{http://www.antoniocorrente.it/wwwsn/images/allegati_articoli/SupplementaryDCM.zip}{all utility matrices}. To get more robust recommendations on the problem at hand, we therefore take into account the 56 multicriteria evaluation combinations using the SMAA methodology and providing the rank acceptability indices and the pairwise winning indices shown in the tables below. Let us remind that the rank acceptability index $b_{k}(a)$ gives the frequency with which the alternative $a$ gets the $k$-th position, with $k=1,\ldots,5$, while the pairwise winning index \citep{LeskinenEtAl2006}, $p(a_i,a_j)$, gives the frequency with which the project $a_{i}$ is preferred to the project $a_{j}$, with $i,j=1,\ldots,5$.

\[
    \left\vert
    \begin{array}{c|ccccc}\hline
                  & b_{1}(\cdot) & b_{2}(\cdot) & b_{3}(\cdot) & b_{4}(\cdot) & b_{5}(\cdot) \\ \hline
        a_{1}     & 0 & 0 & 0 & 0 & 100 \\
        a_{2}     & 30.35 & 58.92 & 10.71 & 0 & 0 \\
        a_{3}     & 0 & 10.71 & 89.28 & 0 & 0 \\
        a_{4}     & 69.64 & 30.35 & 0 & 0 & 0 \\
        a_{5}     & 0 & 0 & 0 & 100 & 0 \\ \hline
        \end{array}
        \right\vert
        \qquad
         \left\vert
    \begin{array}{c|cccccc}\hline
        p(\cdot,\cdot) & a_{1} & a_{2} & a_{3} & a_{4} & a_{5} \\ \hline
        a_{1}     & 0 & 0 & 0 & 0 & 0 \\
        a_{2}     & 100 & 0 & 89.28 & 30.35 & 100 \\
        a_{3}     & 100 & 10.71 & 0 & 0 & 100 \\
        a_{4}     & 100 & 69.64 & 100 & 0 & 100 \\
        a_{5}     & 100 & 0 & 0 & 0 & 0 \\ \hline
        \end{array}
        \right\vert
    \]

Looking at these tables, one can observe that the best project has to be chosen between $a_{4}$ and $a_{2}$ being the only that can take the first ranking position. In particular, $a_{4}$ is a bit in advantage since it has a first rank acceptability index equal to $69.64\%$ against the $30.35\%$ of $a_{2}$. Looking at the worst among the considered projects there is not any doubt about $a_{1}$ since it is always in the last position. In this way, the DM can chose $a_{4}$ or he can wish to investigate a bit more on both projects neglecting all the others.

To conclude this section, let us show how the results change if we admit that the number of blank cards among different levels is not integer as assumed in \cite{CorrenteEtAl2017b}. Let us consider the preference information provided by the DM on criterion $g_{3}$ and represented in the following comparison table: 

 \[
    \left\vert
    \begin{array}{c|ccccccc}\hline
                  & l_{1}(vb) & l_{2}(b) & l_{3}(rb) & l_{4}(a) & l_{5}(rg) & l_{6}(g) & l_{7}(vg) \\ \hline
        l_{1}(vb) &  \blacksquare & \left[1,2\right] & \left[3,4\right] & \left[6,7\right] & \left[9,10\right] & \left[13,14\right] & \left[18,19\right] \\
        l_{2}(b)~ &  & \blacksquare   & \left[1,2\right] & \left[4,5\right] & \left[7,8\right] & \left[11,12\right] & \left[16,17\right] \\
        l_{3}(rb) &  &   &  \blacksquare  & \left[2,3\right] & \left[5,6\right] & \left[9,10\right] & \left[14,15\right] \\
        l_{4}(a)~ &  &   &   & \blacksquare   & \left[2,3\right] & \left[6,7\right] & \left[11,12\right] \\
        l_{5}(rg) &  &   &   &   & \blacksquare   & \left[3,4\right] & \left[8,9\right] \\
        l_{6}(g)~ &  &   &   &   &   & \blacksquare   & \left[4,5\right] \\
        l_{7}(vg) &  &   &   &   &   &    & \blacksquare  \\ \hline
        \end{array}
        \right\vert
    \]

As we already observed in section \ref{Utilities}, solving the $MILP-I$ problem, we find $z^{*}=0$. Therefore, there exists at least one precise comparison table, which elements are integer, compatible with the preferences provided by the DM. Consequently, there exists at least one comparison table  which elements are non-integer, compatible with the same preferences. For example, the following comparison table 

 \[
    \left\vert
    \begin{array}{c|ccccccc}\hline
                  & l_{1}(vb) & l_{2}(b) & l_{3}(rb) & l_{4}(a) & l_{5}(rg) & l_{6}(g) & l_{7}(vg) \\ \hline
        l_{1}(vb) &  \blacksquare & 1.169 & 3.319 & 6.3557 & 9.7306 &  13.7817 & 18.9972 \\
        l_{2}(b)~ &  & \blacksquare   & 1.15 & 4.1866 & 7.5616 & 11.6127 & 16.8282 \\
        l_{3}(rb) &  &   &  \blacksquare  & 2.0367 & 5.4146 & 9.4627 & 14.6782 \\
        l_{4}(a)~ &  &   &   & \blacksquare   & 2.3749 & 6.426 & 11.6415 \\
        l_{5}(rg) &  &   &   &   & \blacksquare   & 3.0511 & 8.2666 \\
        l_{6}(g)~ &  &   &   &   &   & \blacksquare   & 4.2155 \\
        l_{7}(vg) &  &   &   &   &   &    & \blacksquare  \\ \hline
        \end{array}
        \right\vert
    \]
    
\noindent satisfies the consistency Condition \ref{Precies Consistency Condition}, it is concordant with the imprecise preferences provided by the DM and its entries are non-integer. Differently from the case in which the number of blank cards is an integer value and, therefore, the number of comparison tables compatible with the preferences provided by the DM is finite, considering a continuous number of blank cards to be included between two successive levels, there is an infinite number of compatible comparison tables. In particular, defining $P=\{(p,q): p=1,\ldots,6; q=(p+1),\ldots,7\}$, all comparison tables $[e_{pq}]_{(p,q)\in P}$, which elements satisfy these constraints

\begin{subequations}
\begin{gather}
\hspace{1.9cm}{e}_{pk}+{e}_{kq}+1={e}_{pq},  \qquad  (p,k), (k,q), (p,q) \in P, \\
 {e}^{L}_{pq}\leqslant{e}_{pq}\leqslant{e}^{R}_{pq}, \qquad (p,q)\in P, \\
  \hspace{1.5cm}  {e}_{pq}\in \Bbb{R}^{+}_0, \;\;\; \qquad (p,q)  \in P.
\end{gather}
\end{subequations}

 \noindent are compatible with the preference provided by the DM. Therefore, we applied the Hit-And-Run method \citep{Smith1984,TervonenEtAl2012a,VanValkenhoefTervonenPostmus2014}\footnote{See \cite{CorrenteGrecoSlowinski2019} for a detailed description of the application of the HAR method in MCDA.} to sample from the space defined by the constraints above 10,000 precise compatible comparison tables. For each of these comparison tables, we can compute the evaluations for all levels. For example, considering the precise values in the comparison table above and assuming that $u_{3}(vg=7)=100$ and $u_{3}(vb=1)=0$ we get $h=e_{17}+1=19.9972$ and, consequently, $\alpha=\frac{u_{3}(vg=7)-u_{3}(vb=1)}{h}=\frac{100-0}{19.9972}=5.0007$. Therefore,
 
    \begin{itemize}
      \item[~]  $u_3(b=2)= u_{3}(vb=1) + (e_{12}+1)\cdot \alpha= 0 + 2.169\cdot 5.0007= 10.8467$,
      \item[~]  $u_3(rb=3)= u_{3}(vb=1) + (e_{13}+1)\cdot \alpha = 0 + 4.419\cdot 5.0007 = 21.5982$,
      \item[~]  $u_3(a=4)= u_{3}(vb=1) + (e_{14}+1)\cdot \alpha = 0 + 7.3557\cdot 5.0007 = 36.7835$,
      \item[~]  $u_3(rg=5)= u_{3}(vb=1) + (e_{15}+1)\cdot \alpha = 0 + 10.7306\cdot 5.0007 = 53.6605$,
      \item[~]  $u_3(g=6)= u_{3}(vb=1) + (e_{16}+1)\cdot \alpha = 0 + 14.7817\cdot 5.0007 = 73.9188$.
    \end{itemize}

Repeating this procedure for the 10,000 sampled precise comparison tables and considering that on $g_{5}$ there were 8 different comparison tables while for all the other criteria there was just one precise comparison table, we computed the rank acceptability indices and the pairwise winning indices shown in the tables below taking into account all 80,000 multicriteria evaluation combinations:

\[
    \left\vert
    \begin{array}{c|ccccc}\hline
                  & b_{1}(\cdot) & b_{2}(\cdot) & b_{3}(\cdot) & b_{4}(\cdot) & b_{5}(\cdot) \\ \hline
        a_{1}     & 0 & 0 & 0 & 0 & 100 \\
        a_{2}     & 37.45 & 47.44 & 15.10 & 0 & 0 \\
        a_{3}     & 0 & 15.10 & 84.89 & 0 & 0 \\
        a_{4}     & 62.54 & 37.45 & 0 & 0 & 0 \\
        a_{5}     & 0 & 0 & 0 & 100 & 0 \\ \hline
        \end{array}
        \right\vert
        \qquad
         \left\vert
    \begin{array}{c|cccccc}\hline
        p(\cdot,\cdot) & a_{1} & a_{2} & a_{3} & a_{4} & a_{5} \\ \hline
        a_{1}     & 0 & 0 & 0 & 0 & 0 \\
        a_{2}     & 100 & 0 & 84.89 & 37.45 & 100 \\
        a_{3}     & 100 & 15.10 & 0 & 0 & 100 \\
        a_{4}     & 100 & 62.54 & 100 & 0 & 100 \\
        a_{5}     & 100 & 0 & 0 & 0 & 0 \\ \hline
        \end{array}
        \right\vert
    \]

\noindent It is easy to observe that the obtained results are not very far from those computed previously in considering an integer number of blank cards between different levels (for another application of the SMAA approach to the Choquet integral model based on the concept of nonadditive robust ordinal regression see \cite{AngilellaGrecoMatarazzo2010a}).

\section{Conclusions}\label{sec:Conclusions}%%%
%%%%%%%%%%%%%%%%%%%%%%%%%%%%%%%%%%%%%%%%%%%%%%%
\noindent We proposed a new methodology to elicit preference information from a DM to assign values to parameters of preference models in MCDA. These parameters can be weights of criteria in outranking methods, weights of criteria in a SWING based methodology, values assigned by a capacity to subsets of criteria in Choquet integral preference model, marginal value functions. The methodology we presented puts together the simple interpretation and the visual support of the decision deck of cards method and the richer and finer information supplied by comparative judgments representing the difference of appreciation between pairs of elements of the comparison tables approaches. This permits to improve the reliability of the values elicited with our methodology. Taking into account the limited human cognitive capacity we provide also procedures permitting to handle inconsistency of information supplied by the DM as well as incompleteness or imprecise or approximate definition of pairwise comparison tables. We believe that the methodology we are proposing has a great potential because it puts together several basic aspects of very successful approaches in multiple criteria decision aiding, that are:
\begin{itemize}[label={--}]
    \item deck of cards method \citep{FigueiraRoy2002},
    \item pairwise comparison tables from AHP \citep{Saaty1977} and MACBETH \citep{BanaCostaVansnick1994}
    \item SWING method \citep{vonWinterfeldtEd86} to define the weight of criteria in multiattribute value theory approach.
\end{itemize}

The methodology we are proposing is strongly based on the interaction with the DM that can supply the information about values and preferences on which he is convinced enough. Indeed, the pairwise comparison tables we are considering have not to be necessarily complete, which means that the DM is not forced to make comparisons on which she/he is not sure. Moreover, the DM can supply also imprecise information in terms of intervals for the number of cards expressing differences of appreciations of considered elements. The robustness concerns  related to the intervals of number of cards are taken into account in the methodology through the Stochastic Multicriteria Acceptability Analysis (SMAA) \citep{LahdelmaHokkanenSalminen1998}. The methodology we are proposing permits to test the consistency of the information supplied by the DM and, in case of inconsistency, it supplies the DM with explanations of such inconsistency as well as with suggestions to correct it. These are clearly important elements for a discussion between the DM and the analyst that permits the her/him to mature a more comprehensive understanding of the decision problem. 

Due to the above mentioned interesting good properties, we expect therefore that the pairwise comparison table based on the deck of cards approach could be applied to several MCDA decision models. In this article we proposed its application to the Choquet integral decision model permitting to represent interaction between criteria. Another class of decision models to which very naturally our approach can be applied is the family of outranking methods for which pairwise comparison tables based on deck of cards can be used to assess weights of criteria. The methodology we are proposing could also been coupled with the ordinal regression \citep{JacquetLagrezeSiskos1982} and the robust ordinal regression approach \cite{GrecoMousseauSlowinski2008}. Indeed, the pairwise comparison table based on the deck of the cards can be used as a form of preference and value information supplied by the DM in order to assess one value function or a family of compatible value functions.

Together and beyond the methodological development we have just outlined, we expect several applications of the pairwise comparison tables based on the deck of cards in real world decision problems. Indeed, they can support the search of good compromise solutions in situations in which several highly conflicting criteria have to be taken into account such as, for example, in territorial and urban planning, energy system management and sustainable development.

\vspace{0.5cm}

\section*{Acknowledgements}
\addcontentsline{toc}{section}{\numberline{}Acknowledgements}
\noindent Salvatore Corrente and Salvatore Greco gratefully acknowledge the funding by the FIR of the University of Catania BCAEA3, New developments in Multiple Criteria Decision Aiding (MCDA) and their application to territorial competitiveness”, as well as by the research project ``Data analytics for entrepreneurial ecosystems, sustainable development and wellbeing indices" of the Department of Economics and Business of the University of Catania. Salvatore Greco has also benefited of the fund ``Chance'' of the University of Catania. Jos{\'e} Rui Figueira acknowledges the support from the hSNS FCT – Research Project  (PTDC/EGE-OGE/30546/2017), and the FCT grant SFRH/BSAB/139892/2018 under POCH Program.

\addcontentsline{toc}{section}{\numberline{}References}
\bibliographystyle{model2-names}
\bibliography{Bib_Cards,Full_bibliography}

%% not want to use model2-names.bst.



%% \bibitem must have one of the following forms:
%%   \bibitem[Jones et al.(1990)]{key}...
%%   \bibitem[Jones et al.(1990)Jones, Baker, and Williams]{key}...
%%   \bibitem[Jones et al., 1990]{key}...
%%   \bibitem[\protect\citeauthoryear{Jones, Baker, and Williams}{Jones
%%       et al.}{1990}]{key}...
%%   \bibitem[\protect\citeauthoryear{Jones et al.}{1990}]{key}...
%%   \bibitem[\protect\astroncite{Jones et al.}{1990}]{key}...
%%   \bibitem[\protect\citename{Jones et al., }1990]{key}...
%%   \harvarditem[Jones et al.]{Jones, Baker, and Williams}{1990}{key}...
%%


% \end{thebibliography}

\end{document}